\documentclass[a4paper,reqno]{amsart}

\usepackage{amsthm,amssymb,amsmath,mathtools,pdfpages,changes,exscale,cite,color,amsopn,color}
\usepackage[colorlinks=true,pdftex,unicode=true,linktocpage,bookmarksopen,hypertexnames=false]{hyperref}

\newcommand{\free}[1]{\langle#1\rangle}
\newcommand{\tr}{\triangleleft}
\newcommand{\Z}{\mathbb{Z}}
\newcommand{\ov}{\overline}
\newcommand{\GG}{\mathcal{G}}
\renewcommand{\leq}{\leqslant}
\renewcommand{\geq}{\geqslant}

\DeclareMathOperator{\Aut}{Aut}
\DeclareMathOperator{\Fa}{Fa}
\DeclareMathOperator{\gen}{gen}
\DeclareMathOperator{\gr}{gr}
\DeclareMathOperator{\Hol}{Hol}
\DeclareMathOperator{\Img}{Im}
\DeclareMathOperator{\Inj}{Inj}
\DeclareMathOperator{\id}{id}
\DeclareMathOperator{\Ker}{Ker}
\DeclareMathOperator{\Map}{Map}
\DeclareMathOperator{\mpl}{mpl}
\DeclareMathOperator{\op}{op}
\DeclareMathOperator{\Ret}{Ret}
\DeclareMathOperator{\Soc}{Soc}
\DeclareMathOperator{\Sym}{Sym}

\newtheorem{theorem}{Theorem}[section]
\newtheorem{lemma}[theorem]{Lemma}
\newtheorem{proposition}[theorem]{Proposition}
\newtheorem{corollary}[theorem]{Corollary}

\theoremstyle{definition}
\newtheorem{definition}[theorem]{Definition}
\newtheorem{remark}[theorem]{Remark}
\newtheorem{example}[theorem]{Example}

\newtheorem{problem}[theorem]{Problem}
\newtheorem*{convention}{Convention}

\numberwithin{equation}{section}

\subjclass[2010]{16T25, 20F18, 20F19, 20M07}
\keywords{Yang--Baxter equation, set-theoretic solution, nilpotent group, Malcev nilpotent semigroup, multipermutation solution, skew brace}

\begin{document}

\title[On various types of nilpotency of the structure monoid and group]
{On various types of nilpotency of the structure monoid and group of a set-theoretic solution of the Yang--Baxter equation}

\author{F. Ced{\'o} \and E. Jespers \and {\L}. Kubat \and A. Van Antwerpen \and C. Verwimp}

\address[F. Ced\'o]{Departament de Matem\`atiques, Universitat Aut\`onoma de Barcelona, 08193 Bellaterra (Barcelona), Spain}
\email{cedo@mat.uab.cat}

\address[E. Jespers, A. Van Antwerpen, and C. Verwimp]{Department of Mathematics and Data Science,
Vrije Universiteit Brussel, Pleinlaan 2, 1050 Brussel}
\email{eric.jespers@vub.be}
\email{arne.van.antwerpen@vub.be}
\email{charlotte.verwimp@vub.be}

\address[{\L}. Kubat]{Institute of Mathematics, University of Warsaw, Banacha 2, 02-097 Warsaw, Poland}
\email{lukasz.kubat@mimuw.edu.pl}

\date{}

\begin{abstract}
Given a finite bijective non-degenerate set-theoretic solution $(X,r)$ of the Yang--Baxter equation
we characterize when its structure monoid $M(X,r)$ is Malcev nilpotent. Applying this characterization
to solutions coming from racks, we rediscover some results obtained recently by Lebed and Mortier,
and by Lebed and Vendramin on the description of finite abelian racks and quandles.

We also investigate bijective non-degenerate multipermutation (not necessarily finite) solutions
$(X,r)$ and show, for example, that this property is equivalent to the solution associated to the
structure monoid $M(X,r)$ (respectively structure group $G(X,r)$) being a multipermuation solution
and that $G(X,r)$ is solvable of derived length not exceeding the multipermutation level of $(X,r)$
enlarged by one, generalizing results of Gateva-Ivanova and Cameron obtained in the square-free
involutive case. Moreover, we also prove that if $X$ is finite and $G=G(X,r)$ is nilpotent, then the
torsion part of the group $G$ is finite, it coincides with the commutator subgroup $[G,G]_+$ of the
additive structure of the skew left brace $G$ and $G/[G,G]_+$ is a trivial left brace.
\end{abstract}

\maketitle

\section*{Introduction}

Drinfeld \cite{DR1992} suggested to study set-theoretic solutions $(X,r)$ of the Yang--Baxter equation,
that is $X$ is a non-empty set and $r\colon X^2\to X^2$ is a bijective map such that the braid relation
\[r_1\circ r_2\circ r_1=r_2\circ r_1\circ r_2\] holds in $X^{3}$, where $r_1=r \times{\id_X}$ and
$r_2={\id_X}\times r$. We shall write \[r(x,y)=(\lambda_x (y),\rho_y(x))\] for $x,y\in X$. The solution
is said to be involutive if $r^2=\id_{X^2}$. If all maps $\lambda_x$ and $\rho_y$ are bijective then the
solution is said to be 
non-degenerate. 

The description of these solutions has been reduced to a description of associated algebraic objects.
Consequently, set-theoretic solutions of the Yang--Baxter equation are a meeting-ground of mathematical
physics, algebra and combinatorics.

In order to study involutive non-degenerate solutions Rump \cite{Rump2007} introduced left braces and
it has been shown in \cite{BCJ} that all involutive non-degenerate solutions can be obtained from left
braces. In order to deal with non-involutive solutions Guarnieri and Vendramin\cite{GV} introduced skew
left braces and it has been shown in \cite{Ba2018} that all non-degenerate solutions can be obtained from
skew left braces. Furthermore every skew left brace is the permutation group of a non-degenerate solution
\cite{Ba2018}. The structure group \[G=G(X,r)=\gr(X \mid x\circ y=\lambda_x(y) \circ \rho_y(x)\text{ for all }x,y\in X)\]
associated with a non-degenerate solution $(X,r)$ already has been introduced by Etingoff, Schedler and Soloviev
in \cite{ESS} and has since been intensively investigated (see for example \cite{LV,GI2018,GIV,JO2005,gateva2019braided}).
Note that in general the natural map $X\to G(X,r)$ is not injective; if it is, then $(X,r)$ is said to
be an injective solution. In \cite{Yang2018} it has been shown that non-degenerate injective solutions
are isomorphic if and only if their respective structure groups are isomorphic as skew left braces.
In \cite{GIV} Gateva-Ivanova and Van den Bergh introduced the associated structure monoid
\[M=M(X,r)=\langle X \mid x\circ y=\lambda_x(y) \circ \rho_y(x)\text{ for all }x,y\in X\rangle.\]
Clearly $X$ naturally is embedded in $M(X,r)$. In \cite{MGI2008} it is shown that one has an associated
solution $(M,r_M)$ such that $r_M$ restricted to $X^2$ is precisely the original non-degenerate solution
(see also \cite{CJV} for a more general context). The structure monoid thus allows to use algebraic tools
to study non-degenerate solutions. It hence is a fundamental problem to study the algebraic structure
of such monoids. Recall that in \cite{JKV} it has been shown that $M(X,r)$ is a finite module over
an abelian submonoid in case $(X,r)$ is a finite non-degenerate solution.

In this context Ced\'o, Gateva-Ivanova and Smoktunowicz \cite{CGS} proved that if $(X,r)$ is a finite
non-degenerate involutive solution then the structure group $G(X,r)$ is an Engel group (for example
a nilpotent group) if and only if $G(X,r)$ is abelian. Lebed and Mortier \cite{LM} described the finite
quandles with abelian structure group, these are set-theoretic solutions $(X,r)$ with all maps 
$\rho_y=\id_X$ and $\lambda_x(x)=x$ for all $x\in X$ and such that the group
$\gr(\lambda_x:x\in X)\subseteq\Sym(X)$ is abelian.

In \cite{ESS} a special class of non-degenerate solutions has been introduced, the so-called multipermutation
solutions. These subsequently have been intensively studied (see for example
\cite{GIC,BCV,gateva2011quantum,gateva2007set,GI2018}). Jespers and Okni\'nski \cite{JO2005} showed that the 
structure group of a finite involutive non-degenerate multipermutation solution is poly-infinite cyclic and
Bachiller, Ced\'o and Vendramin in \cite{BCV} showed that the converse holds. Ced\'o, Smoktunowicz and Vendramin
\cite{CSV} also have shown that the non-degenerate solution $(B,r_B)$ of the Yang--Baxter equation associated to
a skew left brace $B$ is a multipermutation solution if and only if the skew brace $B$ has finite
multipermutation level. Furthermore, they showed that this is equivalent with the skew left brace is right
nilpotent and the additive structure $(B,+)$ of the skew left brace $(B,+,\circ)$ is nilpotent.

In this paper we continue these investigations. We show two main results on non-degenerate solutions.
The first one is a description of when the structure monoid $M=M(X,r)$ for a finite solution is Malcev
nilpotent. It will be shown that this property is determined by the Malcev nilpotency of cancellative
subsemigroups and a divisibility property on the set $X$ of generators of $M$. In a second main result 
we will show that a non-degenerate solution $(X,r)$ of arbitrary cardinality is a multipermutation solution
if and only if the associated solution $(M,r_M)$ is a multipermutation solution, and this is equivalent with
the associated solution $(G,r_G)$ is a multipermutation solution. Furthermore, if $(X,r)$ is a multipermutation
solution of level $m$ then we prove that the group $G=G(X,r)$ is solvable of derived length bounded by $m+1$.
This extends earlier results of Gateva-Ivanova and Cameron \cite{GIC} on square-free involutive solutions
and of Bachiller, Ced\'o and Vendramin \cite{BCV}. We also show that if $(X,r)$ is a non-degenerate
multipermutation solution on a finite set $X$, then $r$ has even order.

We finish by showing that if $(X,r)$ is a finite multipermutation solution and $G=G(X,r)$ is nilpotent, i.e.,
the additive and multiplicative groups of $G$ are nilpotent and it is right nilpotent as a skew left brace,
then the torsion subgroup $T(G)$ of $(G,\circ)$ is finite and is equal to the additive commutator subgroup $[G,G]_+$ of
the group $(G,+)$, the additive group of the skew left brace $G$, and $\ov{G}=G/[G,G]_+$ is a trivial
left brace. In particular, the image of $(X,r)$ in $(\ov{G},r_{\ov{G}})$ is a trivial solution.

In Section~\ref{sec1} we recall some background and notation on non-degenerate solutions $(X,r)$ of the
Yang--Baxter equation. In particular, we recall that the structure monoid $M(X,r)$ is a regular submonoid
of the holomorph of the structure monoid $A(X,r)=M(X,s)$ of the associated derived solution $(X,s)$ of $(X,r)$.
We also recall the notion for a monoid to be Malcev nilpotent. To deal with our problem we first handle the
case that $M(X,r)=A(X,r)$, i.e., we deal with rack solutions. In Section~\ref{sec3} we then prove necessary
and sufficient conditions for the structure monoid $M(X,r)$ of a finite solution to be Malcev nilpotent.
This will be done by describing a concrete ideal chain that is based on the divisibility properties by
the natural generators $X$ of $M(X,r)$ and $A(X,r)$. In case $(X,r)$ is a multipermutation solution of
level $1$ we give a description of all solutions with Malcev nilpotent structure monoid $M(X,r)$.
In Section~\ref{sec4} we consider multipermutation solutions of arbitrary level. 

\section{Preliminaries}\label{sec1}

In this section we give the necessary background on a bijective non-degenerate set-theoretic solution
of the Yang--Baxter equation and its associated algebraic structures, such as the structure monoid and
structure group as introduced by Etingoff, Schedler and Soloviev \cite{ESS} and Gateva-Ivanova and
Van den Bergh \cite{GIV}. We use the notation of \cite{CJV}. For more details we refer the reader to \cite{CJV,JKV,JKV2}.

Throughout the paper we will use the following notation. For a monoid $S$ and a subset $A$ of $S$ we denote
by $\langle A\rangle$ the submonoid of $S$ generated by $A$. Let $G$ be a group. For every subset $A$ of $G$,
we denote by $\gr(A)$ the subgroup of $G$ generated by $A$. In case $T$ is a semigroup then we denote by
$T^{1}$ the smallest monoid containing $T$. By $\langle X\mid R\rangle$ we denote the monoid presented with
set of generators $X$ and with set of relations $R$. We also use the notation $\gr(X\mid R)$ for the group
presented with set of generators $X$ and with set of relations $R$.

Let $X$ be a non-empty set and let $r\colon X^2\to X^2$ be a map (where we write $r(x,y)=(\lambda_x(y),\rho_y(x))$
for $x,y\in X$). Recall that the pair $(X,r)$ is a bijective non-degenerate (set-theoretic) solution of the
Yang--Baxter equation if $r$ is bijective, all maps $\lambda_x$ and $\rho_y$ are bijective and on
$X^3$ we have
\begin{align}\label{YBE}
	(r\times\id)\circ({\id}\times r)\circ(r\times{\id})=({\id}\times r)\circ(r\times{\id})\circ({\id}\times r).\tag{YBE}
\end{align}
It is well-known and straightforward to check that \eqref{YBE} is equivalent to the following conditions:
\begin{itemize}
	\item[(1)] $\lambda_x\lambda_y=\lambda_{\lambda_x(y)}\lambda_{\rho_y(x)}$,
	\item[(2)] $\lambda_{\rho_{\lambda_x(y)}(z)}\rho_y(x)=\rho_{\lambda_{\rho_x(z)}(y)}\lambda_z(x)$,
	\item[(3)] $\rho_x\rho_y=\rho_{\rho_x(y)}\rho_{\lambda_y(x)}$
\end{itemize}
for all $x,y,z\in X$. In what follows we adopt the following convention.

\begin{convention}
	By a solution of the YBE we will mean a bijective non-degenerate set-theoretic solution of the Yang--Baxter
	equation. For a solution $(X,r)$ of the YBE, we also say that $r\colon X^2\to X^2$ is a solution of the YBE.
\end{convention}

Let $(X,r)$ be a solution of the YBE. The structure monoid $(M,\circ)$ associated to this solution is
\[M=M(X,r)=\langle X\mid x\circ y=\lambda_x(y)\circ\rho_y(x)\text{ for all }x,y\in X\rangle,\]
with identity element denoted by $1$. The structure group $(G,\circ)$ associated to $(X,r)$ is the group
\[G=G(X,r)=\gr( X \mid x\circ y=\lambda_x(y)\circ\rho_y(x) \text{ for all } x,y\in X ).\]
Since $\lambda_x\lambda_y=\lambda_{\lambda_x(y)}\lambda_{\rho_y(x)}$ for all $x,y\in X$, the map
$\lambda\colon X\to\Sym(X)\colon x\mapsto\lambda_x$ extends to a homomorphism $\lambda\colon M\to\Sym(X)$ and induces
a unique homomorphism $\lambda\colon G\to\Sym(X)\colon g\mapsto\lambda_g$, such that
$\lambda_{\iota(x)}(y)=\lambda_x(y)$ for all $x,y\in X$, where $\iota \colon X\to G$ is the natural map.
Similarly, the map $\rho\colon X\to\Sym(X)\colon x\mapsto\rho_x$ extends to an anti-homomorphism $\rho\colon M\to\Sym(X)$
and induces a unique anti-homomorphism $\rho\colon G\to\Sym(X)\colon g\mapsto\rho_g$, such that
$\rho_{\iota(x)}(y)=\rho_x(y)$ for all $x,y\in X$.

It was shown by Gateva-Ivanova and Majid \cite[Theorem~3.6]{MGI2008} that the map $\lambda\colon X\to\Sym(X)$
can be extended to the homomorphism \[\lambda\colon M\to\Map(M,M)\colon m\mapsto\lambda_m,\]
where $\lambda_1=\id_{M}$ and, for $x_1,\dotsc ,x_k,y_1,\dotsc ,y_l\in X$ and $k,l>1$,
\begin{align*}
	\begin{aligned}
		\lambda_{x_1}(1) &=1,\\
		\lambda_{x_1}(y_1\circ\dotsb\circ y_l) &=\lambda_{x_1}(y_1)\circ\lambda_{\rho_{y_1}(x_1)}(y_2\circ\dotsb\circ y_l),\\
		\lambda_{x_1\circ\dots\circ x_k} &=\lambda_{x_1}\circ\dotsb\circ\lambda_{x_k}.
	\end{aligned}\label{lambda1}
\end{align*}
Also, the map $\rho\colon X\to\Sym(X)$ can be extended to the anti-homomorphism
\[\rho\colon M\to\Map(M,M)\colon m\mapsto\rho_m,\]
where $\rho_1=\id_{M}$ and, for $x_1,\dotsc ,x_k,y_1,\dotsc ,y_l\in X$ and $k,l>1$,
\begin{align*}
	\begin{aligned}
		\rho_{x_1}(1) &=1,\\
		\rho_{x_1}(y_1\circ\dotsb\circ y_l) &=\rho_{\lambda_{y_l}(x_1)}(y_1\circ\dotsb\circ y_{l-1})\circ\rho_{x_1}(y_l),\\
		\rho_{x_1\circ\dots\circ x_k} &=\rho_{x_k}\circ\dotsb\circ\rho_{x_1}.
	\end{aligned}\label{erho1}
\end{align*}
Furthermore, \[m_1\circ m_2=\lambda_{m_1}(m_2)\circ\rho_{m_2}(m_1)\] and
\begin{align*}
	\rho_{m_1}(m_2\circ m_3) &=\rho_{\lambda_{m_3}(m_1)}(m_2)\circ\rho_{m_1}(m_3),\\
	\lambda_{m_1}(m_2\circ m_3) &=\lambda_{m_1}(m_2)\circ\lambda_{\rho_{m_2}(m_1)}(m_3)
\end{align*}
for all $m_1,m_2,m_3\in M$. It turns out that $(M,r_M)$, where
\[r_M\colon M\times M\to M\times M\colon (m_1,m_2)\mapsto(\lambda_{m_1}(m_2),\rho_{m_2}(m_1)),\]
is a solution of the YBE (that obviously extends the solution $(X,r)$).

As mentioned in \cite{CJV}, the proof of the above stated result shows that the mappings $\lambda_x$ and $\rho_x$
actually induce left and right actions on $G=G(X,r)$, say
\[\lambda^e\colon G\to\Sym (G)\quad\text{and}\quad\rho^e\colon G\to\Sym (G)\]
such that, furthermore, the mapping
\[r_G \colon G\times G\to G\times G\colon (g_1,g_2)\mapsto(\lambda^e_{g_1}(g_2),\rho^e_{g_2}(g_1))\]
gives a solution $(G,r_G)$ of the YBE. Note that the natural mapping $X\to M(X,r)$ obviously is injective, while,
in general, the natural map $\iota\colon X\to G(X,r)$ is not injective. One says that the solution $(X,r)$ of
the YBE is injective if the map $\iota$ is injective. One obtains that $(G,r_G)$ is an extension of the induced
set-theoretic solution $(\iota(X),r_{\iota(X)})=(\iota(X),r_G|_{\iota(X)^{2}})$, called the injectivization of
$(X,r)$ and sometimes also denoted $\Inj(X,r)$, and $G(X,r)=G(\iota(X),r_{\iota(X)})=G(\Inj(X,r))$.

The left derived solution of $(X,r)$ is the solution $(X,s)$ where
\[s(x,y)=(y,\sigma_y(x))\quad\text{with}\quad\sigma_y(x)=\lambda_y(\rho_{\lambda_x^{-1}(y)}(x)).\]
The structure monoid associated to the left derived solution is called the left derived structure monoid and is defined by
\[A=A(X,r)=\langle X\mid x+\lambda_x(y)=\lambda_x(y)+\lambda_{\lambda_x(y)}(\rho_y(x))\text{ for all }x,y\in X\rangle.\]
Similarly, the right derived structure monoid, which is the structure monoid of the right derived solution $(X,s')$ where
\[s'(x,y)=(\tau_x(y),x)\quad\text{with}\quad\tau_x(y)=\rho_x(\lambda_{\rho_y^{-1}(x)}(y)),\]
is defined as
\[A'=A'(X,r)=\langle X\mid\rho_y(x)\oplus y=\rho_{\rho_y(x)}(\lambda_x(y))\oplus\rho_y(x)\text{ for all }x,y\in X\rangle.\]
We denote the identity element of $A(X,r)$ (respectively of $A'(X,r)$) by $0$ (respectively by $0'$).
The map $\lambda \colon X\to\Sym(X)$ can be extended to the homomorphism
\[\lambda'\colon M\to\Aut(A,+)\colon m\mapsto\lambda'_m,\] where $\lambda'_1=\id_A$, $\lambda'_m(0)=0$
for all $m\in M$, and
\[\lambda'_{ x_1\circ\dotsb\circ x_k}(y_1+\dotsb+y_l)=
\lambda_{x_1}\dotsm\lambda_{x_k}(y_1)+\dotsb+\lambda_{x_1}\dotsm\lambda_{x_k}(y_l)\]
for all $x_1,\dotsc,x_k,y_1,\dotsc,y_l\in X$ and $k,l\geq 1$. Also, the map $\rho\colon X\to\Sym(X)$
can be extended to the anti-homomorphism \[\rho'\colon M\to\Aut (A',\oplus)\colon m\mapsto\rho'_m,\]
where $\rho'_1=\id_{A'}$, $\rho'_m(0')=0'$ for all $m\in M$, and
\[\rho'_{ x_1\circ\dotsb\circ x_k}(y_1\oplus\dotsb\oplus y_l)
=\rho_{x_k}\dotsm\rho_{x_1}(y_1)\oplus\dots\oplus\rho_{x_k}\dotsm\rho_{x_1}(y_l)\]
for all $x_1, \dotsc, x_k, y_1, \dotsc, y_l\in X$ and $k,l\geq 1$.

In \cite{JKV,JKV2} it has been shown that for a solution $(X,r)$ of the YBE, the structure monoid $M(X,r)$
and the associated derived structure monoids are strongly linked (this has been extended in \cite{CJV}
to arbitrary solutions). Indeed,
\begin{itemize}
	\item[(i)] there is a unique bijective $1$-cocycle $\pi\colon M(X,r)\to A(X,r)$ with respect to the
	left action $\lambda'$ such that $\pi(x)=x$ for all $x\in X$. So, for all $m_1,m_2\in M$ we have
	$\pi(m_1\circ m_2)=\pi(m_1)+\lambda'_{m_1}(\pi(m_2))$.
	\item[(ii)] there is a unique bijective $1$-cocycle $\pi'\colon M(X,r)\to A'(X,r)$ with respect to the
	right action $\rho'$ such that $\pi'(x)=x$ for all $x\in X$. So, for all $m_1,m_2\in M$ we have
	$\pi'(m_2\circ m_1)=\rho'_{m_1}(\pi'(m_2))\oplus\pi'(m_1)$.
\end{itemize}
Furthermore, the mapping \[f\colon M(X,r)\to A(X,r)\rtimes\Img(\lambda')\colon m\mapsto (\pi (m),\lambda_m')\]
is a monoid monomorphism and the mapping \[f'\colon M(X,r)\to A'(X,r)^{\op}\rtimes\Img(\rho')\colon m\mapsto (\pi'(m),\rho'_m)\]
is a monoid anti-monomorphism, where $A'(X,r)^{\op}$ is the opposite monoid of $A'(X,r)$. Hence, for any $m_1,m_2\in M$,
\begin{align*}
	(\pi(m_1),\lambda'_{m_1})(\pi(m_2),\lambda'_{m_2}) &=(\pi(m_1)+\lambda'_{m_1}(\pi(m_2)),\lambda'_{m_1}\lambda'_{m_2}),\\
	(\pi'(m_1),\rho'_{m_1})(\pi'(m_2),\rho'_{m_2}) &=(\rho'_{m_1}(\pi'(m_2))\oplus\pi'(m_1),\rho'_{m_1}\rho'_{m_2}).
\end{align*}

The structure group associated to the left derived solution is called the left derived structure group and is defined by
\[A_{\gr}(X,r)=\gr( X \mid x+\lambda_x(y)=\lambda_x(y)+\lambda_{\lambda_x(y)}(\rho_y(x))\text{ for all }x,y\in X ).\]
Notice that $A_{\gr}(X,r)=G(X,s)$. Similarly, the right derived structure group, which is the structure group of the
right derived solution $(X,s')$, is defined by
\[A_{\gr}'(X,r)=\gr(X\mid\rho_y(x)\oplus y=\rho_{\rho_y(x)}(\lambda_x(y))\oplus\rho_y(x)\text{ for all }x,y\in X ).\]
As above, we have $A_{\gr}'(X,r)=G(X,s')$. Again there exists a bijective $1$-cocycle $\pi \colon G(X,r)\to A_{\gr}(X,r)$
(for simplicity we also denote this map again by $\pi$) such that $\pi(x)=x$ for all $x\in X$, and a morphism of groups
\[\lambda^{e'}\colon G(X,r)\to\Aut(A_{\gr}(X,r),+)\colon g\mapsto\lambda^{e'}_g,\]
such that the mapping \[G(X,r)\to A_{\gr}(X,r)\rtimes\Img(\lambda^{e'})\colon g\mapsto(\pi (g),\lambda^{e'}_g)\]
is a group monomorphism. The bijective $1$-cocycle transfers the additive group structure on $A_{\gr}(X,r)$ to an
additive structure on $G=G(X,r)$, which we also denote by $+$. Hence $G$ is equipped with two operations $+$ and
$\circ$ that are related by the following ``distributive'' property: \[a\circ(b+c)=(a\circ b)-a+(a\circ c)\] for all
$a,b,c\in G$. So $(G,+,\circ)$ is a skew left brace, as introduced by Guarnieri and Vendramin in \cite{GV}.
In case $(G,+)$ is abelian this is simply called a left brace, a notion introduced by Rump in \cite{Rump2007}.

The following proposition clarifies the link between the maps $\lambda' \colon M\to\Aut (A,+)$ and
$\lambda\colon M\to\Map(M,M)$ as well as between the maps $\rho'\colon M\to\Aut(A',\oplus)$ and $\rho\colon M\to\Map(M,M)$.

\begin{proposition}\label{prop:rel}
	Let $(X,r)$ be a solution of the YBE. Then
	\[\pi(\lambda_{m_1}(m_2))=\lambda'_{m_1}(\pi(m_2))\quad\text{and}\quad\pi'(\rho_{m_1}(m_2))=\rho'_{m_1}(\pi'(m_2))\]
	for all $m_1,m_2\in M(X,r)$.
	\begin{proof}
		We prove the first part only; the second part can be proven similarly. To do so, we first prove that
		$\pi(\lambda_x(m))=\lambda'_x(\pi(m))$ for all $x\in X$ and $m\in M(X,r)$, by induction on the length
		$|m|$ of $m$. If $|m|=1$ then \[\pi(\lambda_x(m))=\lambda_x(m)=\lambda_x(\pi(m))=\lambda'_x(\pi(m)).\]
		Suppose we have proven the result for words in $M(X,r)$ of length at most $k$. Let $m\in M(X,r)$ be an
		element of length $k+1$. Write $m=y\circ m'$, where $y\in X$ and $m'\in M(X,r)$ with $|m'|=k$. Then
		\begin{align*}
			\pi(\lambda_x(m)) &=\pi(\lambda_x(y\circ m'))\\
			& =\pi(\lambda_x(y)\circ\lambda_{\rho_y(x)}(m'))\\
			& =\pi(\lambda_x(y))+\lambda'_{\lambda_x(y)}(\pi(\lambda_{\rho_y(x)}(m')))\\
			& =\lambda'_x(\pi(y))+\lambda'_{\lambda_x(y)}(\lambda'_{\rho_y(x)}(\pi(m')))\\
			& =\lambda'_x(\pi(y))+\lambda'_x(\lambda'_y(\pi(m')))\\
			& =\lambda'_x(\pi(y)+\lambda'_y(\pi(m')))\\
			& =\lambda'_x(\pi(y\circ m'))\\
			& =\lambda'_x(\pi(m)).
		\end{align*}
		Using that both $\lambda$ and $\lambda'$ are homomorphisms, we obtain $\pi(\lambda_{m_1}(m_2))=\lambda'_{m_1}(\pi(m_2))$
		for all $m_1,m_2\in M(X,r)$, and the result follows.
	\end{proof}
\end{proposition}

Note that the group $\gr(\lambda_x:x\in X)\subseteq\Sym(X)$ is isomorphic to the group
$\gr(\lambda_m:m\in M)=\gr(\lambda_x:x\in X)\subseteq\Sym(M)$. Indeed, via the $1$-cocycle $\pi$, we may identify $M$
with $A=A(X,r)$, and by Proposition~\ref{prop:rel}, $\lambda_m$ is then identified with $\lambda'_m\in\Aut(A,+)$. Since
$X$ generates $A$, the map $\lambda_x\in \Sym(X)$ determines $\lambda'_x\in\Aut(A,+)$, and thus also determines
$\lambda_x\in \Sym(M)$, by its identification with $\lambda'_x$. Now the map
\[\{\lambda_x\in \Sym(X):x\in X\}\to\{\lambda_x\in \Sym(M):x\in X\},\]
defined by $\lambda_x\mapsto\lambda_x$, induces an isomorphism of the groups
$\gr(\lambda_x:x\in X)$ and $\gr(\lambda_m:m\in M)$. Similarly one can see that the group
$\gr(\rho_x:x\in X)\subseteq\Sym(X)$ is isomorphic to the group $\gr(\rho_m:m\in M)=\gr(\rho_x:x\in X)\subseteq\Sym(M)$.

In \cite{JKV,JKV2} it is proven that if $(X,r)$ is a finite solution of the YBE then $M=M(X,r)$ is a finite
(left and right) module over an abelian normal submonoid $T$ of $M$ (which may be embedded in $A(X,r)$), i.e.,
$M=\bigcup_{f\in F} Tf=\bigcup_{f\in F}fT$ for some finite subset $F$ of $M$. Hence, $G(X,r)$ is (finitely generated)
abelian-by-finite.

Note that for every solution $(X,r)$ of the YBE, $(X,r^{-1})$ also is a solution of the YBE. Write
$r^{-1}(x,y)=(\hat{\lambda}_x(y),\hat{\rho}_y(x))$ for $x,y\in X$. We define some types of permutation
groups associated to the solution $(X,r)$:
\begin{align*}
	\GG_{\gen}(X,r) &=\gr((\lambda_x,\rho_x^{-1},\hat{\lambda}_x,\hat{\rho}_x^{-1}):x\in X)\subseteq\Sym(X)^4,\\
	\GG_{\lambda,\rho}(X,r) &=\gr((\lambda_x,\rho_x^{-1}):x\in X)\subseteq\Sym(X)^2,\\
	\GG_{\lambda,\hat{\lambda}}(X,r) &=\gr((\lambda_x, \hat{\lambda}_x):x\in X)\subseteq\Sym(X)^2,\\
	\GG_{\lambda}(X,r) &=\gr(\lambda_x:x\in X)\subseteq\Sym(X),\\
	\GG_{\rho}(X,r) &=\gr(\rho_x:x\in X)\subseteq\Sym(X).
\end{align*}
Bachiller in \cite[Definition 3.10]{Ba2018} (or \cite[Definition 2.1.13]{DBThesis}) defined the following permutation
group: \[\gr((\lambda_x, \tilde{g_x}^{-1}):x\in X)=\{ (\lambda_a, \tilde{g}_a^{-1}):a\in G(X,r)\}\subseteq\Sym(X)^2,\]
where $\tilde{g}$ is defined as follows: \[\tilde{g}_a(y)=\rho_{(\lambda^e_y)^{-1}(a)}(y)\] for all $a\in G(X,r)$
and $y\in X$. Note that $\hat{\lambda}_{\lambda_y(x)}(\rho_x(y))=y$ for all $x,y\in X$. Hence
$\hat{\lambda}_z(\rho_{\lambda_y^{-1}(z)}(y))=y$ and thus $\hat{\lambda}_z^{-1}(y)=\rho_{\lambda_y^{-1}(z)}(y)=\tilde{g}_z(y)$
for all $y,z\in X$. Thus the permutation group in the sense of Bachiller is in our notation the group $\GG_{\lambda,\hat{\lambda}}(X,r)$.

\begin{lemma}\label{permutationlemma}
	Let $(X,r)$ be a solution of the YBE. Then
 	\begin{alignat*}{2}
		\lambda_a^{-1}(x) &=\hat{\rho}_{(\hat{\lambda}^e_x)^{-1}(a)}(x),\quad
		& \hat{\lambda}_a^{-1}(x) &=\rho_{(\lambda^e_x)^{-1}(a)}(x),\\
		\rho_a^{-1}(x) &=\hat{\lambda}_{(\hat{\rho}^e_x)^{-1}(a)}(x),\quad
		& \hat{\rho}_a^{-1}(x) &=\lambda_{(\rho^e_x)^{-1}(a)}(x)
 	\end{alignat*}
 	for all $a\in G(X,r)$ and $x\in X$.
 	\begin{proof}
		By \cite[Lemma~2.1.12]{DBThesis}, the map $\tilde{g}\colon G(X,r)\to\Sym(X)$, defined by $\tilde{g}(a)=\tilde{g}_a$
		and $\tilde{g}_a(x)=\rho_{(\lambda^e_x)^{-1}(a)}(x)$ for all $a\in G(X,r)$ and $x\in X$, is an anti-homomorphism
		of groups. Similarly, one verifies that the map $\tilde{f}\colon G(X,r)\to\Sym(X)$, defined by $\tilde{f}(a)=\tilde{f}_a$,
		where $\tilde{f}_a(x)=\lambda_{(\rho^e_x)^{-1}(a)}(x)$ for all $a\in G(X,r)$ and $x\in X$, is a homomorphism of groups.
		Note that the map $\hat{\lambda}^{-1}\colon G(X,r)\to\Sym(X)\colon a\mapsto \hat{\lambda}_a^{-1}$ is an anti-homomorphism
		of groups. Since \[\hat{\lambda}_x^{-1}(y)=\rho_{\lambda_y^{-1}(x)}(y)=\tilde{g}_x(y)\] for all $x,y\in X$, we have that
		$\hat{\lambda}^{-1}=\tilde{g}$, i.e., \[\hat{\lambda}_a^{-1}(x)=\rho_{(\lambda^e_x)^{-1}(a)}(x)\] for all $a\in G(X,r)$
		and $x\in X$. We also have that the map $\hat{\rho}^{-1}\colon G(X,r)\to\Sym(X)\colon a\mapsto \hat{\rho}_a^{-1}$ is a homomorphism
		of groups. Since \[\hat{\rho}_x^{-1}(y)=\lambda_{\rho_y^{-1}(x)}(y)=\tilde{f}_x(y)\] for all $x,y\in X$, we have that
		$\hat{\rho}^{-1}=\tilde{f}$, i.e., \[\hat{\rho}_a^{-1}(x)=\lambda_{(\rho^e_x)^{-1}(a)}(x)\] for all $a\in G(X,r)$ and
		$x\in X$. This proves two of the equalities in the statement of the result. The other two equalities follow similarly. 
 	\end{proof}
\end{lemma}

\begin{lemma}\label{Permutation}
	Let $(X,r)$ be a solution of the YBE. Then the groups $\GG_{\gen}(X,r)$, $\GG_{\lambda,\rho}(X,r)$
	and $\GG_{\lambda,\hat{\lambda}}(X,r)$ are isomorphic.
	\begin{proof}
		Note that the maps \[h_1\colon G(X,r)\to\GG_{\lambda, \hat{\lambda}}(X,r)\colon a\mapsto (\lambda_a, \hat{\lambda}_a),\]
		\[h_2\colon G(X,r)\to\GG_{\lambda, \rho}(X,r)\colon a\mapsto (\lambda_a, \rho_a^{-1})\]
		and \[h\colon G(X,r)\to\GG_{\gen}(X,r)\colon a\mapsto (\lambda_a,\rho_a^{-1}, \hat{\lambda}_a,\hat{\rho}_a^{-1})\]
		are epimorphisms of groups.
		
		First we shall see that $\rho^e_x(\Ker(\lambda))=\Ker(\lambda)$ for all $x\in X$. To do so, let us fix
		$a\in\Ker(\lambda)$ and $x\in X$. We have
		\[\lambda_x=\lambda_a\lambda_x=\lambda_{\lambda^e_a(x)}\lambda_{\rho^e_x(a)}=\lambda_x\lambda_{\rho^e_x(a)}\]
		and \[\lambda_{x^{-1}}=\lambda_a\lambda_{x^{-1}}=\lambda_{\lambda^e_a(x^{-1})}\lambda_{\rho^e_{x^{-1}}(a)}=
		\lambda_{x^{-1}}\lambda_{\rho^e_{x^{-1}}(a)}.\] Hence $\lambda_{\rho^e_x(a)}=\id$ and
		$\lambda_{(\rho^e_x)^{-1}(a)}=\id$, and thus $\rho^e_x(\Ker(\lambda))=\Ker(\lambda)$.
		
		Similarly one proves that
		\[\lambda^e_x(\Ker(\rho))=\Ker(\rho),\quad \hat{\rho}^e_x(\Ker(\hat{\lambda}))=\Ker(\hat{\lambda}),\quad
		\hat{\lambda}^e_x(\Ker(\hat{\rho}))=\Ker(\hat{\rho})\] for all $x\in X$.
		
		Now, we shall see that \[\Ker(h_1)=\Ker(\lambda)\cap\Ker(\rho)=\Ker(\rho)\cap\Ker(\hat{\rho}),\] and therefore
		$\Ker(h_1)=\Ker(h_2)=\Ker(h)$. Indeed, let $a\in\Ker(h_1)=\Ker(\lambda)\cap\Ker(\hat{\lambda})$. For every $x\in X$,
		we have $\hat{\lambda}_{(\hat{\rho}_x^e)^{-1}(a)}=\id$ and $ \lambda_{(\rho_x^e)^{-1}(a)}=\id$. Thus by
		Lemma~\ref{permutationlemma}, \[\rho_a^{-1} (x)=\hat{\lambda}_{(\hat{\rho}^e_x)^{-1}(a)}(x)=x\quad\text{and}\quad
		\hat{\rho}_a^{-1}(x)=\lambda_{(\rho^e_x)^{-1}(a)}(x)=x\] for all $x\in X$. This shows that $\Ker(h_1)
		\subseteq\Ker(\rho)\cap\Ker(\hat{\rho})$. The other inclusion follows by a symmetric argument. Hence $\Ker(h_1)=\Ker(h)$.
		Thus, we obtain that \[\GG_{\lambda,\hat{\lambda}}{(X,r)}\cong G(X,r)/\Ker(h_1)=G(X,r)/\Ker(h)\cong\GG_{\gen}(X,r).\]
		
		Note that we have also proven that $\Ker(h_1)\subseteq\Ker(\rho)\cap\Ker(\lambda)=\Ker(h_2)$. Let $b\in\Ker(h_2)$.
		As $\Ker(\rho)$ is $\lambda_x^e$-invariant, for every $x\in X$, we have that $\rho_{(\lambda^e_x)^{-1}(b)}=\id$.
		Therefore $\hat{\lambda}_b^{-1}(x)=\rho_{(\lambda^e_x)^{-1}(b)}(x)=x$ for all $x\in X$. This shows that
		$\Ker(h_2)\subseteq\Ker(h_1)$, and thus $\Ker(h_1)=\Ker(h_2)$. Hence
		\[\GG_{\lambda,\rho}(X,r)\cong G(X,r)/\Ker(h_2)=G(X,r)/\Ker(h)\cong\GG_{\gen}(X,r),\] and the result follows.
	\end{proof}
\end{lemma}

\begin{remark}\label{Kernel}
	Note that in the proof of Lemma~\ref{Permutation} we have shown that
	$\Ker(\lambda)\cap\Ker(\rho)=\Ker(\hat{\lambda})\cap\Ker(\hat{\rho})$.
	Indeed, we find that \[\Ker(\lambda)\cap\Ker(\rho)=\Ker(\lambda)\cap\Ker(\rho)\cap\Ker(\hat{\lambda})\cap\Ker(\hat{\rho}),\]
	which by symmetry between $(X,r)$ and $(X,r^{-1})$ entails the result.
\end{remark}

\begin{definition}
	Let $(X,r)$ be a solution of the YBE. We define the permutation group $\GG(X,r)$
	of $(X,r)$ as \[\GG(X,r)=\GG_{\lambda,\rho}(X,r).\]
\end{definition}

Note that in case the solution $(X,r)$ is involutive, we have $\GG(X,r)=\GG_{\lambda,\rho}(X,r)\cong \GG_{\lambda}(X,r)=\GG_{\rho}(X,r)$.

\begin{remark}\label{rem:h2}
	Let $(X,r)$ be a solution of the YBE. Put $G=G(X,r)$. We know that $(G,+,\circ)$ is a skew left brace.
	Since \[\lambda_a^e(b)=-a+a\circ b\quad\text{and}\quad \rho_a^e(b)=(b^{-1}+a)^{-1}\circ a\] for $a,b\in G$,
	it follows that the socle of $G$ is
	\begin{align*}
		\Soc(G) & =\{a\in G:a\circ b=a+b=b+a\text{ for all }b\in G\}\\
		& =\{ a\in G:\lambda^e_a=\id\text{ and }\rho^e_a=\id\}.
	\end{align*}
	The socle is an ideal of the skew left brace $G$. Both its additive and multiplicative groups are abelian;
	in particular it is a left brace. Note that both $\lambda_a^e$ and $\rho_a^e$ are maps $G\to G$. If
	$(X,r)$ is an injective solution (for example if $r$ is an involution), i.e., the natural map $\iota\colon X\to G$
	is injective, then $\lambda_a^e=\id_G$ if and only if $\lambda_a=\lambda_a^e|_X=\id_X$. In this case $\Soc(G)=\Ker(h_{2})$
	and thus $G/\Soc(G)\cong\GG_{\lambda,\rho}(X,r)$. However, in general, we know that for $a\in G$, $\lambda_a=\id_X$ implies
	$\lambda_a^e=\id_G$; and similarly for the $\rho$-maps. Hence, \[\Ker(h_2) \subseteq \Soc(G)\] and thus $G/\Soc(G)$ is an
	epimorphic image of \[G/\Ker(h_2)\cong\GG_{\lambda,\rho}(X,r).\] In \cite[Example~3.12]{Ba2018}
	(or \cite[Example~2.1.15]{DBThesis}) Bachiller gave an example where $\GG_{\lambda,\hat{\lambda}}(X,r)$ and 
	$G/\Soc(G)$ are not isomorphic. In fact, in this example
	\[G/\Soc(G)\cong\Z/2\Z\quad\text{and}\quad\GG(X,r)\cong\GG_{\lambda,\hat{\lambda}}(X,r)\cong\GG_{\lambda,\rho}(X,r)\cong\Z.\]
	Another example, with the same idea, is the following. Let $n$ be a positive integer and $X=\Z/n\Z$. Let $r\colon X^2\to X^2$
	be defined by $r(x,y)=(y+1,x+1)$ for all $x,y\in X$. In $G=G(X,r)$ we have $x\circ (x-1)=x\circ (x+1)$ for each $x\in X$ and
	thus $\iota(x-1)=\iota(x+1)$ for all $x\in X$, where $\iota\colon X\to G$ is the natural map. If $n$ is odd, then $G$ is free
	abelian of rank $1$, and $\Soc(G)=G$. Hence $G/\Soc(G)=0$ if $n$ is odd. If $n$ is even then $G=\gr(0,1\mid 0\circ 0=1\circ 1)$
	and $\Soc(G)=\gr(0\circ 0, 1\circ 0, 0\circ 1)$. Hence $G/\Soc(G)\cong\Z/2\Z$ if $n$ is even. In both cases, $\GG(X,r)\cong \Z/n\Z$.
	
	Since $\Ker(h_2)$ is a normal subgroup of $G(X,r)$, it is easy to see that $\Ker(h_2)$ is an ideal of the skew left
	brace $G(X,r)$. This allows to define an addition on $\GG(X,r)$ by
	$(\lambda_a,\rho_a^{-1})+(\lambda_b,\rho_b^{-1})=(\lambda_{a+b},\rho_{a+b}^{-1})$ for all $a,b\in G(X,r)$.
	Then $(\GG(X,r),+,\circ)$ is a skew left brace (see \cite[Theorem 3.11]{Ba2018} or \cite[Theorem 2.1.14]{DBThesis}).
\end{remark}

\section{Malcev nilpotency of $A(X,r)$} \label{sec2}

We begin with stating a known result on the structure of nilpotent structure groups.
The second part is due to Ced\'o, Gateva-Ivanova and Smoktunowicz \cite{CGS} and Lebed and Vendramin \cite{LV}.

For a group $G$ we denote by $T(G)$ the set consisting of the elements of finite order (also called the torsion elements).

\begin{lemma}\label{lemmaabelian}
	Let $(X,r)$ be a finite solution of the YBE and $G=G(X,r)$.
	\begin{enumerate}
		\item If the group $G$ is nilpotent then $G$ is finite-by-(free abelian). In particular, in this case,
		it is a finite conjugacy group (i.e., $G$ has finite commutator subgroup).
		\item If $G$ is torsion-free then $G$ is nilpotent if and only if $G$ is abelian, or equivalently
		the injectivization $\Inj(X,r)$ of $(X,r)$ is the trivial solution. Hence, if $(X,r)$ is a finite
		non-degenerate involutive solution then $G$ is nilpotent if and only if $G$ is abelian.
	\end{enumerate}
	\begin{proof}
		(1) We know that $G$ is abelian-by-finite and finitely generated. Assume $G$ also is nilpotent. Then $T=T(G)$
		is a finite characteristic subgroup of $G$. Since $G/T$ is torsion-free and abelian-by-finite, it follows that $G/T$
		is finitely generated abelian. Hence a free abelian group. (2) Since $G\cong G(\Inj (X,r))$ it is sufficient
		to prove this for injective solutions $(X,r)$. It was shown in \cite{JKVV} that $G$ is torsion-free if and
		only if $\Inj(X,r)$ is involutive. The result now follows from the involutive case.
	\end{proof}
\end{lemma}

Malcev showed that nilpotency of groups can be defined via some specific identities, the so-called Malcev identities.
For convenience of the reader we recall the definition of a Malcev nilpotent semigroup \cite{Malcev}. Let $F$ denote
the free semigroup on $\{x,y,z_n:n\geq 1\}$. For non-negative integers $n$ define Malcev's words
\[x_n=x_n(x,y;z_1,\dotsc,z_n)\in F\quad\text{and}\quad y_n=y_n(x,y;z_1,\dotsc,z_n)\in F\] recursively as
\begin{alignat*}{2}
	x_0 &=x,\quad & y_0 &=y,\\
	x_{n+1} &=x_n z_{n+1} y_n,\quad & y_{n+1} &=y_n z_{n+1} x_n.
\end{alignat*}
A semigroup $S$ is said to be Malcev nilpotent of nilpotency class $n$ (or, simply, nilpotent of class $n$)
if $n$ is the smallest non-negative integer such that \[x_n(s,t;u_1,\dotsc,u_n)=y_n(s,t;u_1,\dotsc,u_n)\]
in $S$ for all $s,t\in S$ and $u_1,\dotsc,u_n\in S^1$. Recall that a group $H$ is Malcev nilpotent of class
$n$ if and only if $H$ is nilpotent (in ordinary sense) of class $n$.

Let $(X,r)$ be a solution of the YBE. Note that if $M=M(X,r)$ is nilpotent of class $n$ then, as an epimorphic image,
the group $\gr(\lambda'_m:m\in M)$ also is nilpotent and thus so is the isomorphic group $\gr(\lambda_x:x\in X)$.

In this section we determine when this structure monoid $M(X,r)$ is nilpotent in case $M(X,r)=A(X,r)$ or $M(X,r)=A'(X,r)$,
that is, in case all $\lambda_x=\id$ or in case all $\rho_x=\id$. The characterization is given for the latter case, i.e.,
for $r(x,y)=(\lambda_x (y),x)$. For such solutions it is more customary the write $\lambda_x(y)$ as $y\tr x$ and $r$ as
$r_{\tr}$, and to use the context of racks. Recall (see for example \cite{AndGra2003}) that a rack is a set $X$ with a
binary operation denoted $\tr$ satisfying for any $x,y,z\in X$:
\begin{enumerate}
	\item[(R1)] $(x\tr y)\tr z=(x\tr z)\tr(y\tr z)$ for all $x,y,z\in X$ (right self-distributivity),
	\item[(R2)] the right translation $X\to X\colon x\mapsto x\tr y$ is a bijection for each $y\in X$.
\end{enumerate}
If, furthermore,
\begin{enumerate}
	\item[(R3)] $x\tr x=x$ for each $x\in X$ (idempotence)
\end{enumerate}
then $(X,\tr)$ is called a quandle (see for example \cite{Joyce1982}).

Note that if $(X,r)$ is a solution of the YBE such that $\rho_x=\id$ for all $x\in X$, then $(X,\tr)$, where
$y\tr x=\lambda_x(y)$, is a rack. Conversely, if $(X,\tr)$ is a rack, then $(X,r_{\tr})$, where $r_{\tr}(x,y)=(y\tr x,x)$,
is a solution of the YBE.

Let $(X,\tr)$ be a rack. Consider its associated solution $(X,r_{\tr})$ and write $\lambda_x(y)=y\tr x$. Thus the
associated permutation group is $\GG(X,r_{\tr})\cong\GG_{\lambda}(X,r_{\tr})$, and, for simplicity, we identify
these groups for this class of solutions, that is, $\GG(X,r_{\tr})=\gr(\lambda_x:x\in X)$. We also denote
$\GG(X,r_{\tr})$, $M(X,r_{\tr})$ and $G(X,r_{\tr})$ by $\GG(X,\tr)$, $M(X,\tr)$ and $G(X,\tr)$, respectively.

\begin{proposition}\label{racknilpotent}
	Assume that $(X,\tr)$ is a rack. If the permutation group $\GG(X,\tr)$ is nilpotent of class $n$
	then the structure monoid $M(X,\tr)$ of the solution $(X,r_\tr)$ of the YBE associated to $(X,\tr)$
	is Malcev nilpotent of class not exceeding $n+2$. Similarly, the structure group $G(X,\tr)$ of
	$(X,r_\tr)$ is nilpotent of class not exceeding $n+2$.
	\begin{proof}
		Let $M=M(X,\tr)$. Assume that $\GG(X,\tr)$ is nilpotent of class $n$. Choose $a,b,a_1,\dotsc,a_{n+2}\in M$
		and define
		\begin{align*}
			x &=x_{n}(a,b;a_1,\dotsc,a_{n})\in M,\\
			y &=y_{n}(a,b;a_1,\dotsc,a_{n})\in M.
		\end{align*}
		Since $\GG(X,\tr)$ is Malcev nilpotent of class $ n$, we get that $\lambda_x=\lambda_y$.
		Since $r_\tr$ is a solution of the YBE, we get $\lambda_y\lambda_x=\lambda_x\lambda_y=\lambda_{\lambda_x(y)}\lambda_x$.
		Hence $\lambda_{\lambda_x(y)}=\lambda_y$. Now, let $z=\lambda_y(y)\in M$. Because
		$\lambda_y\lambda_y=\lambda_{\lambda_y(y)}\lambda_y=\lambda_z\lambda_y$, we have $\lambda_y=\lambda_z$.
		Thus $\lambda_z\lambda_x=\lambda_x\lambda_z=\lambda_{\lambda_x(z)}\lambda_x$ and, in consequence,
		$\lambda_{\lambda_x(z)}=\lambda_z$. Furthermore,
		\[x\circ z=\lambda_x(z)\circ x=\lambda_{\lambda_x(z)}(x)\circ\lambda_x(z)=\lambda_z(x)\circ\lambda_x(z),\]
		which leads to
		\begin{align*}
			y\circ x\circ x\circ y &=\lambda_y(x)\circ y\circ\lambda_x(y)\circ x
			=\lambda_y(x)\circ\lambda_y(\lambda_x(y))\circ y\circ x\\
			&=\lambda_z(x)\circ\lambda_x(z)\circ y\circ x
			=x\circ z\circ y\circ x=x\circ \lambda_y(y)\circ y\circ x\\
			&=x\circ y\circ y\circ x.
		\end{align*}
		Then, using the previous equality, we obtain
		\begin{align*}
			\begin{aligned}
				y\circ x\circ a\circ x\circ y &=y\circ \lambda_x(a)\circ x\circ x\circ y
				=\lambda_y(\lambda_x(a))\circ y\circ x\circ x\circ y\\
				&=\lambda_x(\lambda_y(a))\circ x\circ y\circ y\circ x
				=x\circ \lambda_y(a)\circ y\circ y\circ x\\
				&=x\circ y\circ a\circ y\circ x
			\end{aligned}\label{abeliannilpoteneq}
		\end{align*}
		for all $a\in M$. Finally, the last equality leads to
		\begin{align*}
			y_{n+2}(a,b;a_1,\dotsc,a_{n+2}) &=y\circ (a_{n+1}\circ x)\circ a_{n+2}\circ (x\circ a_{n+1})\circ y\\
			&=y\circ x\circ \lambda_x^{-1}(a_{n+1})\circ a_{n+2}\circ \lambda_x(a_{n+1})\circ x\circ y\\
			&=(x\circ y)\circ \lambda_x^{-1}(a_{n+1})\circ a_{n+2}\circ \lambda_x(a_{n+1})\circ (y\circ x)\\
			&=x\circ \lambda_y(\lambda_x^{-1}(a_{n+1}))\circ y\circ a_{n+2}\circ y\circ \lambda_y^{-1}(\lambda_x(a_{n+1}))\circ x\\
			&=x\circ a_{n+1}\circ y\circ a_{n+2}\circ y\circ a_{n+1}\circ x\\
			&=x_{n+2}(a,b;a_1,\dotsc,a_{n+2}).
		\end{align*}
		So $M$ is Malcev nilpotent of class at most $ n+2$. A similar argument shows that $G(X,\tr)$ is Malcev nilpotent
		of class $\leq n+2$, and thus it is nilpotent of class $\leq n+2$.
	\end{proof}
\end{proposition}

It is worth to add that Proposition~\ref{racknilpotent} may be strengthened in case $(X,\tr)$ is a quandle.
Moreover, in this case, it is possible to provide a simpler proof.

\begin{corollary}
	Let $(X,\tr)$ be a quandle. If the permutation group $\GG(X,\tr)$ is nilpotent of class $n$ then
	the structure monoid $M(X,\tr)$ and structure group $G(X,\tr)$ of the solution $(X,r_\tr)$
	are Malcev nilpotent of class at most $n+1$.
	\begin{proof}
		Let $M=M(X,\tr)$. Choose $a,b,a_1,\dotsc,a_{n+1}\in M$ and define
		\begin{align*}
			x_i &=x_i(a,b;a_1,\dotsc,a_i)\in M,\\
			y_i &=y_i(a,b;a_1,\dotsc,a_i)\in M
		\end{align*}
		for $1\leq i\leq n+1$. Since $\GG(X,\tr)$ is Malcev nilpotent of class $n$, we obtain that
		$\lambda_{x_n}=\lambda_{y_n}$. Moreover, the fact that $(X,\tr)$ is a quandle assures that
		its associated solution $(X,r_\tr)$ is square-free,  which implies that $\lambda_m(m)=m$
		for each $m\in M$. Taking the above into account we get
		\[x_n\circ y_n=\lambda_{x_n}(y_n)\circ x_n=\lambda_{y_n}(y_n)\circ x_n=y_n\circ x_n\]
		in $M$ and consequently
		\begin{align*}
			x_{n+1} &=x_n\circ a_{n+1}\circ y_n=\lambda_{x_n}(a_{n+1})\circ x_n\circ y_n\\
			&=\lambda_{y_n}(a_{n+1})\circ y_n\circ x_n=y_n\circ a_{n+1}\circ x_n=y_{n+1}.\qedhere
		\end{align*}
	\end{proof}
\end{corollary}

In a recent paper Lebed and Mortier \cite{LM} describe all finite quandles $(X,\tr )$ with abelian
structure group. In particular, these quandles are abelian, i.e., $(a\tr b)\tr c=(a \tr c)\tr b$
for all $a,b,c\in X$. Equivalently the associated permutation group $\GG(X,\tr )$ is abelian. Furthermore,
the structure group $G=G(X,\tr)$ of an abelian quandle is presented as a central extension of a free
abelian group by an explicit finite abelian group. The latter easily can be proven as follows. For this
note that in $G$ we have $\lambda_x(y)=xyx^{-1}$ for all $x,y\in X$. The map $\lambda\colon X\to\Sym(X)\colon x\mapsto\lambda_x$
induces a unique homomorphism \[\lambda\colon G\to\GG(X,\tr)\colon g\mapsto\lambda_g.\]
Note that, for $x_1,\dots ,x_n\in X$ and $\varepsilon_1,\dotsc,\varepsilon_n\in\{-1,1\}$,
$x_1^{\varepsilon_1}\dotsm x_n^{\varepsilon_n}\in\Ker(\lambda)$ if and only if
$\lambda^{\varepsilon_1}_{x_1}\dotsm\lambda^{\varepsilon_n}_{x_n}=\id$. Since
\[x_1^{\varepsilon_1}\circ\dotsb\circ x_n^{\varepsilon_n}\circ y
=\lambda_{x_1}^{\varepsilon_1}\dotsm\lambda_{x_n}^{\varepsilon_n}(y)\circ
x_1^{\varepsilon_1}\circ\dotsb\circ x_n^{\varepsilon_n}\] for all $y\in X$,
it follows that $\Ker(\lambda)\subseteq Z(G)$. Hence the group $G/Z(G)$ is a homomorphic image of the group
$G/\Ker(\lambda)\cong\GG(X,\tr)$. However, it may happen that $\Ker(\lambda)\ne Z(G)$ (compare with Remark~\ref{rem:h2}).
To provide an explicit example of this phenomenon consider $X=\{1,2\}$ and define $x\tr y=\sigma(x)$ for $x,y\in X$,
where $\sigma$ is the unique non-trivial permutation of $X$. Then
\[G=G(X,\tr)=\gr(1,2\mid 1\circ 1=2\circ 1=2\circ 2=1\circ 2)\cong\Z.\]
Hence $G=Z(G)$, but $G/\Ker(\lambda)\cong\GG(X,\tr)\cong\Z/2\Z$. Thus $\Ker(\lambda)\ne Z(G)$.

So we have shown the following result (compare with Proposition~\ref{racknilpotent}).

\begin{corollary}\label{racknilpotentgroup}
	Assume that $(X,\tr)$ is a rack and $G=G(X,\tr)$. Then the group $G/Z(G)$ is a homomorphic
	image of the group $\GG=\GG(X,\tr)$. In particular, $G$ is nilpotent if and only if $\GG$ is nilpotent,
	and in this case the class of nilpotency of $G$ is equal to  or exceeds by one  the class of nilpotency of $\GG$.
	Furthermore, $G$ is solvable if and only if $\GG$ is solvable, and in this case  the derived length of $G$
	is equal to or exceeds by one the derived length of $\GG$.
\end{corollary}

\begin{corollary}[{Lebed and Mortier \cite[Theorem 3.2]{LM}}]\label{lebedmortier}
	Assume $(X,\tr)$ is a finite abelian rack and $G=G(X,\tr)$. Then $G$ is a finite conjugacy
	group with periodic subgroup $T(G)=G'$ and $G/G'$ is a free abelian group of rank at most $|\iota (X)|$.
\end{corollary}

The first part follows from the following proposition.

\begin{proposition}\label{pro:alltorsionaxr}
	Let $(X,r)$ be a finite solution of the YBE. If $A=A'_{\gr}(X,r)$ then $T(A)$, the set of torsion elements
	of $A$, coincides with the finite group $[A,A]$. Also $T(A_{\gr}(X,r))=[A_{\gr}(X,r),A_{\gr}(X,r)]$.
	\begin{proof}
		We only prove the first part of the result.
		
		Because $A$ is a finitely generated finite conjugacy group, the commutator subgroup $[A,A]$ of $A$
		is a finite group (see for example \cite{Neumann}).
		
		Recall that $A$ is the structure group of the right derived solution $(X,s')$ of $(X,r)$,
		where $s'(x,y)=(\tau_x(y),x)$ and $\tau_x(y)=\rho_x(\lambda_{\rho_y^{-1}(x)}(y))$ for all $x,y\in X$.
		Let $\tau\colon A\to\Sym(X)\colon a\mapsto \tau_a$ be the unique homomorphism such that $\tau_{\iota(x)}=\tau_x$
		for all $x\in X$, where $\iota\colon X\to A$ is the natural map. Consider the equivalence relation $\approx$
		on $X$, where $x \approx y$ for $x,y\in X$ if there exists $a\in A$ such that $\tau_a(x)=y$. Write
		$[x]\in\ov{X}=X/{\approx}$ for the $\approx$-class of $x\in X$. Consider the free abelian group $F=\Fa(\ov{X})$
		on $\ov{X}$. Note that $F$ is the structure group of the solution $(\ov{X},\ov{s}')$ of the YBE,
		where $\ov{s}'([x],[y])=([\tau_x(y)],[x])=([y],[x])$ for all $x,y\in X$. Because the map $(X,s')\to(\ov{X},\ov{s}')$,
		defined as $x\mapsto [x]$ for $x\in X$, is an epimorphism of solutions, there is a unique morphism of groups
		$\varphi\colon A\to F$ such that $\varphi(x)=[x]$ for all $x\in X$. Clearly $\varphi$ factors uniquely through
		a homomorphism $\ov{\varphi}\colon A/[A,A]\to F$. On the other hand, $x\oplus y=\tau_x(y)\oplus x$ in $A$
		and thus $\tau_x(y)\ominus y\in [A,A]$ for all $x,y\in X$. Hence the map $\ov{X}\to A/[A,A]$, defined by
		$[x]\mapsto \ov{x}$ for $x\in X$, is well defined. Hence there exists a unique homomorphism
		$\psi\colon F\to A/[A,A]$ such that $\psi([x])=\ov{x}$ for all $x\in X$. Clearly $\psi$ is
		the inverse of $\ov{\varphi}$. Therefore $T(A)=[A,A]$ and the result is proved.
	\end{proof}
\end{proposition}

We now consider the natural action of the permutation group $\GG(X,\tr)$ on the set $X$. Let $X=X_1\sqcup\dotsb\sqcup X_r$
be a decomposition of $X$ into orbits with respect to the action of $\GG(X,\tr)$ on $X$. So, for $1\leq i\leq r$, we have
$x,y\in X_i$ if and only if there exists $g\in \GG(X,\tr)$ such that $g(x)=y$.

\begin{lemma}\label{lemma:q3}
	Let $(X,\tr)$ be an abelian rack. If $x,y\in X$ belong to the same $\GG(X,\tr)$-orbit then $\lambda_x=\lambda_y$.
	\begin{proof}
		Because the rack $(X,\tr)$ is abelian we have \[(x \tr y) \tr z=(x \tr z) \tr y=(x \tr y) \tr (z \tr y)\]
		for all $x,y,z\in X$. Hence, $\lambda_z(\lambda_y(x))=\lambda_{\lambda_y(z)}(\lambda_y(x))$ for all $x,y,z\in X$,
		and thus $\lambda_z\lambda_y=\lambda_{\lambda_y(z)}\lambda_y$ for all $y,z\in X$. As $\lambda_y$ is bijective,
		we get $\lambda_z=\lambda_{\lambda_y(z)}$, and the result follows.
	\end{proof}
\end{lemma}

\begin{lemma}\label{lemma:q4}
	Let $(X,\tr)$ be a finite abelian rack with $\GG(X,\tr)$-orbits $X_1,\dotsc, X_r$. If $1\leq i,j\leq r$
	and $x_i\in X_i$, then the map $\lambda_{x_i}|_{X_j}$ is a permutation of $X_j$ consisting of disjoint
	cycles all of the same length.
	\begin{proof}
		Let $\sigma=\lambda_{x_i}|_{X_j}=\sigma_1\dotsb \sigma_s$ be the decomposition of the permutation
		$\sigma$ of $X_j$ as a product of disjoint cycles. We may assume that $\sigma \neq \id$ and that
		$\sigma_1$ has minimal length, say $n$. This implies that $\sigma^n=\sigma_1^n\dotsb\sigma_s^n$ has
		a fixed point, say $x_j\in X_j$. We claim that $\sigma^n$ is the identity map. Indeed, take $x\in X_j$.
		As $x$ is in the same orbit as $x_j$, there exists $g\in \GG(X,\tr)$ such that $g(x_j)=x$.
		Because $\GG(X,\tr)$ is abelian, this implies $\sigma^n(x)=\sigma^n(g(x_j))=g(\sigma^n(x_j))=g(x_j)=x$.
		Thus all the disjoint cycles of $\sigma$ must have length $n$.
	\end{proof}
\end{lemma}

With the assumptions as in Lemma~\ref{lemma:q4}, consider the following subgroups of $G=G(X,\tr)$:
\[G_{i}=\gr(X_i) \subseteq G.\] Note that in $G_i$ we have $x\circ y=\lambda_x(y)\circ x$
for all $x,y\in X_i$. We claim that $G_i$ is an abelian group.

\begin{lemma}\label{lemma:q5}
	Let $(X,\tr)$ be an abelian finite rack. Then the groups $G_1,\dotsc,G_r$ are abelian.	
	\begin{proof}
		Note that in $G=G(X,\tr)$ we have $x\circ x=\lambda_x(x)\circ x$ for all $x\in X$. Hence,
		$x=\lambda_x(x)$ in $G$. Fix $1\leq i\leq r$. It is enough to prove that all generators
		of $G_i$ commute. So, let $x,y\in X_i$. Because of Lemma~\ref{lemma:q3},
		we get $x\circ y=\lambda_x(y)\circ x=\lambda_y(y)\circ x=y\circ x$, as desired.
	\end{proof}
\end{lemma}

Let $1 \leq i,j \leq r$. For any $x\in X_i$ and $y\in X_j$ define the commutator
\[g_{x,y}=[x,y]=x\circ y\circ x^{-1}\circ y^{-1}\in G(X,\tr).\] As $g_{x,y}=\lambda_x(y)\circ y^{-1}\in G_j$
and $g_{x,y}=x\circ \lambda_y(x)^{-1}\in G_i$, we obtain $g_{x,y}\in G_i \cap G_j$. By Lemma~\ref{lemma:q5},
the groups $G_i$ and $G_j$ are abelian. So, $g_{x,y}$ is central in both $G_i$ and in $G_j$. We claim that
$g_{x,y}$ is also central in $G(X,\tr)$.

\begin{lemma}\label{lemma:q6}
	With the notation as above $g_{x,y}$ is central in $G=G(X,\tr)$ and $g_{x,y}=g_{x',y'}$
	if $x,x'\in X_i$ and $y,y'\in X_j$. We simply denote $g_{x,y}$ as $g_{ij}$.
	\begin{proof}
		First we shall show that $g_{x,y}=g_{x',y'}$ for $x,x'\in X_i$ and $y,y'\in X_j$. Since
		\[g_{x,y}=x\circ y\circ x^{-1}\circ y^{-1}=\lambda_x(y)\circ y^{-1}=x\circ \lambda_y(x)^{-1},\]
		we get $g_{x,y}=\lambda_x(y)\circ y^{-1}=\lambda_{x'}(y)\circ y^{-1}=g_{x',y}$ and similarly
		$g_{x',y}=x'\circ \lambda_y(x')^{-1}=x'\circ \lambda_{y'}(x')^{-1}=g_{x',y'}$. Thus $g_{x,y}=g_{x',y'}$.
		
		Using this observation we shall prove that $g_{x,y}$ commutes with each generator $z\in X$. Indeed, we have
		\begin{align*}
			z\circ g_{x,y} &=z\circ [x,y]\circ z^{-1}\circ z=[z\circ x\circ z^{-1},z\circ y\circ z^{-1}]\circ z\\
			&=[\lambda_z(x),\lambda_z(y)]\circ z=g_{\lambda_z(x),\lambda_z(y)}\circ z=g_{x,y}\circ z.\qedhere
		\end{align*}
	\end{proof}
\end{lemma}

Because of Corollary~\ref{lebedmortier} we obtain the following result.

\begin{corollary}[{Lebed and Vendramin \cite[Theorem 8.15]{LV}}]
	Let $(X,\tr)$ be an abelian finite rack. If $G=G(X,\tr)$ then \[G'=\gr(g_{ij}:1\leq i,j\leq r)=T(G)\subseteq Z(G).\]
	In particular, $G$ is abelian if and only if $G$ is free abelian, or equivalently, $G$ is a torsion-free group.
\end{corollary}

This result has been proven in general in \cite{JKV,JKV2} for arbitrary finite bijective non-degenerate solutions
$(X,r)$: the monoid $M(X,r)$ is free abelian if and only if the algebra $KM(X,r)$ over an arbitrary field $K$
is a domain if and only if the monoid $M(X,r)$ is cancellative. It is easy to see that in this case we get that
$KG(X,r)$ is a domain. Hence, from the positive solution of the zero divisor problem for polycyclic-by-finite groups,
it follows that the latter is equivalent with $G(X,r)$ being a torsion-free group.

\begin{proposition}\label{prop:rackdesc}
	Finite abelian racks on a set $X$ are in one-to-one correspondence with partitions $X=X_1\sqcup\dotsb\sqcup X_r$
	of $X$ together with families of permutations $f_{ij}\in \Sym(X_i)$ for $1 \leq i,j \leq r$ such that:
	\begin{enumerate}
		\item $f_{ij}f_{ik}=f_{ik}f_{ij}$ for all $1\leq i,j,k\leq r$,
		\item if $\GG_i=\gr(f_{ij}:1 \leq j \leq r)$ then $\GG_ix_i=X_i$ for each $1\leq i\leq r$ and $x_i\in X_i$
		(here $\GG_ix_i$ denotes the orbit of $x_i$ with respect to the action of $\GG_i$ on $X_i$),
		\item if $g\in\GG_i$ for some $1\leq i\leq r$ has a fixed point then $g=\id$.
	\end{enumerate}
	Moreover, the decomposition and permutations above correspond to an abelian quandle provided $f_{ii}=\id$ for each $1\leq i\leq r$.
	\begin{proof}
		The result is a consequence of what we have already shown. If $(X,\tr)$ is an abelian rack then we have a decomposition
		$X=X_1\sqcup\dotsb\sqcup X_r$ of $X$ into orbits with respect to the action of $\GG(X,\tr)$ on $X$. In particular, each
		$\lambda_x$ for $x\in X$ preserves the components of this decomposition, that is $\lambda_x(X_i)=X_i$ for $1\leq i\leq r$.
		Therefore, if $1\leq i,j\leq r$ then choosing $x_j\in X_j$ we may define $f_{ij}\in\Sym(X_i)$ as $f_{ij}=\lambda_{x_j}|_{X_i}$;
		Lemma~\ref{lemma:q3} shows that the permutation $f_{ij}$ is well-defined, that is it does not depend on the representative $x_j$
		of the orbit $X_j$. Since $\lambda_{x_j}\lambda_{x_k}=\lambda_{x_k}\lambda_{x_j}$ for all $x_j\in X_j$ and $x_k\in X_k$, we get
		$f_{ij}f_{ik}=f_{ik}f_{ij}$ for all $1\leq i,j,k\leq r$. Moreover, if $x_i\in X_i$ then
		\[X_i=\GG(X,\tr)x_i=\{(g|_{X_i})(x_i):g\in\GG(X,\tr)\}=\GG_ix_i.\]
		Further, if $g\in \GG_i$ satisfies $g(x_i)=x_i$ and $x\in X_i$ then writing $x=f(x_i)$ for some $f\in \GG_i$, we obtain 
		$g(x)=g(f(x_i))=f(g(x_i))=f(x_i)=x$, as desired. Finally, if $(X,\tr)$ is a quandle then $\lambda_x(x)=x$ for $x\in X_i$
		yields $f_{ii}(x_i)=x_i$ and thus, by the previous, $f_{ii}(x)=x$ for each $x\in X_i$, that is $f_{ii}=\id$.
		
		Conversely, having a decomposition $X=X_1\sqcup\dotsb\sqcup X_r$ of $X$ and a family of permutations $f_{ij}\in\Sym(X_i)$
		for $1\leq i,j\leq r$ satisfying conditions (1)--(3), we may define an abelian rack structure on $X$ by declaring that
		$x\tr y=f_{ij}(x)$ for $x,y\in X$ provided $x\in X_i$ and $y\in X_j$. Indeed, if $x\in X_i$, $y\in X_j$ and $z\in X_k$
		for some $1\leq i,j,k\leq r$ then
		\begin{align*}
			(x\tr y)\tr z &=f_{ij}(x)\tr z=f_{ik}(f_{ij}(x))=f_{ij}(f_{ik}(x))\\
			&=f_{ik}(x)\tr f_{jk}(y)=(x\tr z)\tr(y\tr z).
		\end{align*}
		Moreover, if $y\in X_j$ then the map $f\colon X\to X$, defined as $f(x)=x\tr y$, is bijective because $f(x)=f_{ij}(x)$
		for $x\in X_i$ (that is $f|_{X_i}=f_{ij}$). Finally, if $f_{ii}=\id$ for each $1\leq i\leq r$ then $(X,\tr)$ is a quandle,
		because if $x\in X_i$ then $x\tr x=f_{ii}(x)=x$.
		
		It is also easy to check that the correspondence between abelian rack structures on $X$ and decompositions of $X$ together
		with families of maps $f_{ij}$ satisfying all conditions stated in proposition is in fact a one-to-one correspondence.
		Thus the result is proved.
	\end{proof}
\end{proposition}

In \cite[Theorem 2.3]{LM} Lebed and Mortier have obtained a combinatorial description of families of permutations
satisfying the requirements in Proposition~\ref{prop:rackdesc}; and thus they obtained a full description
of all finite abelian racks. This combinatorial description is in terms of $r$-tuples of lower-triangular matrices
with non-negative entries. Quandles corresponding to such $r$-tuples are called in \cite{LM} the filtered-permutation quandles.

\begin{problem}
	Describe the finite quandles $(X,\tr)$ with permutation group $\GG(X,\tr)$
	a nilpotent group of class $2$, or more general, a metabelian group.
\end{problem}

A natural problem is to investigate arbitrary finite solutions $(X,r)$ of the YBE with $\GG(X,r)$ an abelian group.

\section{Malcev nilpotency of $M(X,r)$}\label{sec3}

Let $(X,r)$ be a finite solution of the YBE (recall that, by our convention, we mean that the solution $(X,r)$
is non-degenerate and bijective). As mentioned earlier, the structure monoid is a finite module over an abelian
submonoid and thus $M(X,r)$ is a linear monoid, i.e., a submonoid of the multiplicative monoid of a matrix ring
over a field. Okni{\'n}ski \cite{OknLinear} (and later Jespers and Riley \cite{JR2006}) gave a criterion for a
finitely generated linear semigroup $S$ to be nilpotent. This criterion is based on information of certain ideal
chains of $S$  (actually ideal chains with factors that are either power nilpotent or uniform subsemigroups of
completely $0$-simple inverse semigroups) of which the existence follows from the fact that $M(X,r)$ (or $KM(X,r)$,
where $K$ is a field) is Noetherian. In order to describe when the structure monoid $M(X,r)$ is a Malcev nilpotent
monoid we will give a very concrete description of such an ideal chain. It hence will also give an independent
proof of the previous for the structure monoids $M(X,r)$. This is what we first deal with in this section.

Put $X=\{ x_1, \dotsc, x_n\}$. So $M=M(X,r)=\langle x_1,\dotsc,x_n \rangle $ and $A=A(X,r)=\langle a_1,\dotsc,a_n\rangle$,
where $a_i=\pi (x_i)$ for $1\leq i \leq n$. Further we have a monoid embedding
$f\colon M\to A(X,r)\rtimes\Img(\lambda')\colon m\mapsto (\pi (m), \lambda_{m}')$.
Abusing notation, we will identify $m$ with $f(m)$, i.e., we will write $m=(\pi(m),\lambda_m')$.
For simplicity reasons, for $a\in A$, we will write $\lambda_a'$ for $\lambda'_{\pi^{-1}(a)}$.
So we simply may write \[M=\{(a,\lambda_a'):a\in A\}=\free{x_i=(a_i,\lambda_{a_i}'):1\leq i \leq n}.\]
Therefore, we have a mapping \[\lambda'\colon A\to\Aut(A,+)\colon a\mapsto\lambda_a',\]
and, for $a,b\in A$, \[\lambda_{a+\lambda_a'(b)}'=\lambda_a'\circ\lambda_b' .\]

For a subset $B$ of $A$ we put, as in \cite{JKV}, $B^e=\{(b,\lambda_b'):b\in B\}$.

We will construct an ideal chain in $M $ based on divisibility elements of $X$, the generators of $M$.
This idea has been used in earlier work on monoids of $I$-type and also on monoids of skew and quadratic
type (see for example \cite{JespersOkninksiBook,JO2006,GJO2003,JVC2017}). It also has been used in \cite{JKV,JKV2}
to determine the prime ideals of $M$ and of its algebra $KM$. However, in order to determine when $M$ is nilpotent,
we need to get more detail on the ideal chain constructed from the divisibility by generators.

Recall that an element $s$ in a monoid $S$ is left divisible by $t\in S$ if $s=tt'$ for some $t'\in S$.
Similarly one defines right divisibility. If all elements of $S$ are normalizing (i.e., $Ss=sS$ for all $s\in S$),
such as in the monoid $A$, then left and right divisibility are the same; in this case we simply use the terminology divisible
and we write $t\mid s$. Now, note that in $M$ an element $(a,\lambda_a')$ is left divisible by a generator $x_i=(a_i,\lambda_{a_i}')$
if and only if $a$ is divisible by $a_i$. So, left divisibility in $M$ by elements of $X$ can be transferred to divisibility in $A$
by elements of $\{ a_1, \dotsc, a_n\}$, the generators of $A$.

For $1\leq i \leq n$ put
\begin{align*}
	M_i=\{(a,\lambda'_a)\in M :{} & (a,\lambda_a') \text{ is left divisible by at least } i \text{ different}\\
	{} & \text{generators amongst } x_1,\dotsc, x_n\}
\end{align*}
Clearly, \[M_i=A_i^e=\{ (a,\lambda_a'):a\in A_i\}\] with
\begin{align*}
	A_i=\{a\in A : {} & a \text{ is divisible by at least $i$ different}\\
	{} & \text{generators amongst }a_1,\dotsc,a_n\}.
\end{align*}
Note $a\in A$ being divisible by $a_i$ in $A$ means $a=a_i+ b$ for some $b\in A$, or equivalently
\[(a,\lambda_a')=(a_i,\lambda_{a_i}')((\lambda_{a_i}')^{-1}(b),(\lambda_{a_i}')^{-1}\lambda_a')
=x_i ((\lambda_{a_i}')^{-1}(b),(\lambda_{a_i}')^{-1} \lambda_a').\]
Although $(a,\lambda_a')$ is left divisible by $x_i$ this does not mean that $(a,\lambda_a')$ is right divisible by $x_i$.

As stated in \cite{JKV} (it is easy to verify that) each $M_i$ is a two-sided ideal of $M$. Hence we get in $M$ the ideal chain
\begin{equation}
	\varnothing=M_{n+1}\subseteq M_n\subseteq M_{n-1}\subseteq \dotsb \subseteq M_1\subseteq M_0=M.\tag{$*$}\label{Chain1}
\end{equation}

Next we will refine the above chain \eqref{Chain1}. We will show that there exists ideals $B_i,U_i$ of $M$ satisfying
\[M_{i+1}\subseteq B_i \subseteq U_i\subseteq M_i\] and such that
\begin{enumerate}
	\item $B_i/M_{i+1}$ and $M_i/U_i$ are power nilpotent semigroups (if $M_i/M_{i+1}$ is power nilpotent then we take $B_i=U_i=M_i$),
	\item if $M_i/M_{i+1}$ is not power nilpotent then $U_i\setminus B_i$ is a disjoint union of semigroups $S_1,\dotsc,S_m$ such that 
	$S_kS_l\subseteq M_{i+1}$ for $k\ne l$,
	\item and each $(S_i\cup M_{i+1})/M_{i+1}$ is a uniform subsemigroup of a completely $0$-simple inverse semigroup.
\end{enumerate}
Recall that a completely $0$-simple inverse semigroup is a semigroup of the form $\mathcal{M}^{0}(C,r,r,I)$,
where $C$ is a group and $I$ is the $r\times r$ identity matrix, i.e., this is the semigroup of all $r\times r$
matrices with entries in the group with zero $C\cup\{\theta\}$ (here $\theta$ denotes the zero element) which have
at most one non-zero entry. A subsemigroup $T$ of $\mathcal{M}^{0}(C,r,r, I)$ is said to be uniform if each
$\mathcal{H}$-class (i.e., all the matrices with non-zero entries in a fixed $(i,j)$ spot) of $\mathcal{M}^{0}(C,r,r, I)$
intersects non-trivially $T$ and the maximal subgroups of $\mathcal{M}^{0}(C,r,r, I)$ are generated by their intersection
with $T$. We make the agreement that some ideals in the chain can be empty.

Fix $i$ with $1\leq i \leq n$. We need to introduce some notations. Let
\[\mathcal{L}=\{ Y \subseteq \{a_1,\dotsc,a_n\}:|Y|=i\}.\]
For $Y,Z\in \mathcal{L}$ put \[M_{YZ}=A_{YZ}^e\] with
\[A_{YZ}=\{a\in A\setminus A_{i+1}:y\mid a\text{ for all }y\in Y\text{ and }\lambda_a'(Z)=Y\}.\]
Notice that if $a\in A_{YZ}$ then $a\in\free{Y}$ and $z\nmid a$ if $z\in \{a_1,\dotsc,a_n\}\setminus Y$.
Also note that some elements of $A_i \cap\free{Y}$ might belong to $A_{i+1}$.

For $Y\in\mathcal{L}$ put \[M_{Y*}=\bigcup_{Z\in\mathcal{L}}M_{YZ}\quad\text{and}\quad M_{*Y}=\bigcup_{Z\in\mathcal{L}}M_{ZY}.\]

\begin{lemma}\label{lemma1}
	The following properties hold for $Y\in \mathcal{L}$.
	\begin{enumerate}
		\item $M_{Y*}\cup M_{i+1}$ is a right ideal of $M$.
		\item $M_{*Y}\cup M_{i+1}$ is a left ideal of $M$.
 	\end{enumerate}
 	\begin{proof}
		(1) Let $Z\in \mathcal{L}$ and $(a,\lambda_a')\in M_{YZ}$. Then, for $(b,\lambda_b')\in M$, we have
		$(a,\lambda_a')(b,\lambda_b')=(a+\lambda_a' (b),\lambda_a' \lambda_b')$. If $a+\lambda_a'(b)\notin\free{Y}$,
		i.e., $a+\lambda_a'(b)$ is divisible by some $z\in \{a_1,\dotsc,a_n\}\setminus Y$ then $a+\lambda_a'(b)\in A_{i+1}$
		and thus $(a,\lambda_a')(b,\lambda_b')\in M_{i+1}$. Otherwise, $a+\lambda_a'(b)\in A_{i}\setminus A_{i+1}$.
		Furthermore, \[(\lambda'_{a+\lambda'_a(b)})^{-1}(Y)=(\lambda'_b)^{-1}(\lambda'_a)^{-1}(Y)=(\lambda'_b)^{-1}(Z).\]
		Thus $(a,\lambda_a') (b,\lambda_b')\in M_{Y(\lambda'_b)^{-1}(Z)}\subseteq M_{Y*}$, as desired.
		
		(2) Let $(a,\lambda_a')\in M_{ZY}$, where $Z\in \mathcal{L}$, and let $(b,\lambda_b')\in M$. Then
		$(b,\lambda_b')(a,\lambda_a')=(b+\lambda_b'(a), \lambda_b' \lambda_a')$. If $b+\lambda_b'(a)\in A_{i+1}$
		then $(b,\lambda_b') (a,\lambda_a')\in M_{i+1}$. So, suppose $b+\lambda_b'(a)\in A_{i}\setminus A_{i+1}$
		and thus $(b,\lambda_b') (a,\lambda_a')\in M_{i} \setminus M_{i+1}$. Clearly, $\lambda_b'(a)$ is divisible
		by all elements of $\lambda_b'(Z)$. Furthermore, we have
		$(\lambda_b'\lambda_a')^{-1}(\lambda_b'(Z))=(\lambda'_a)^{-1}(\lambda'_b)^{-1}(\lambda_b'(Z))=(\lambda'_a)^{-1}(Z)=Y$.
		Hence, $(b,\lambda_b') (a,\lambda_a')\in M_{\lambda_b'(Z)Y}\subseteq M_{*Y}$, as desired.
	\end{proof}
\end{lemma}

From \cite{JKV} we know that there exists $d\geq 2$ such that $da\in Z(A)$, the center of $A$, and $\lambda_{da}=\id$ for all
$a\in A$. For $Y\in \mathcal{L}$ we put \[a_Y=\sum_{y\in Y} dy\in A\quad\text{and}\quad m_Y=(a_Y,\lambda_{a_Y}')=(a_Y,\id)\in M.\]
Note that $a_Y$ is divisible by all elements of $Y$ but it could be divisible by more than $i$ generators, i.e.,
$a_Y$ could belong to $A_{i+1}$. 
Also, 
\begin{equation}\label{ia_Y}
	ka_Y = \sum_{y\in Y} kdy\quad\text{and}\quad m_Y^k=(ka_Y,\id),
\end{equation}
for any positive integer $k$.

\begin{lemma}\label{lemma2}
	Let $Y\in \mathcal{L}$. If $a_Y\in A_{i+1} $ then following properties hold:
	\begin{enumerate}
		\item $(M_{YY}\cup M_{i+1})/M_{i+1}$ is nil.
		\item $(M_{*Y}\cup M_{i+1})/M_{i+1}$ is a nil left ideal of $M/M_{i+1}$.
		\item $(M_{Y*}\cup M_{i+1})/M_{i+1}$ is a nil right ideal of $M/M_{i+1}$.
	\end{enumerate}
	Hence, $B_i=M_{i+1}\cup\bigcup_{Y:a_Y\in A_{i+1}}(M_{*Y}\cup M_{Y*})$ is in the radical of $M/M_{i+1}$.
	\begin{proof}
		(1) Let $(a,\lambda_a')\in M_{YY}$. Then, for any positive integer $k$,
		\[(a,\lambda_a')^{k}=(a+\lambda'_a(a)+\dotsb+(\lambda'_a)^{k-1}(a),(\lambda'_a)^k).\]
		Because, $\lambda'_a (Y)=Y$ we get that $(a,\lambda'_a)^{k}\in M_{YY}\cup M_{i+1}$.
		Since each $(\lambda'_a)^i(a)$ is divisible by all elements of $Y$ and because each element of $A$
		is normalizing, we get that $a+\lambda'_a(a)+\dotsm +(\lambda'_a)^{k-1}(a)$ is divisible by $a_Y$ for
		a large enough $k$. Hence, it then follows that $(a,\lambda'_a)^{k}\in m_YM\subseteq M_{i+1}$. Therefore,
		$(M_{YY}\cup M_{i+1})/M_{i+1}$ is nil.
		
		(2) Let $Z\in \mathcal{L}$. Assume that $Z\ne Y$ and $(a,\lambda'_a)\in M_{ZY}$. Then
		$(a,\lambda'_a )(a,\lambda'_a)=(a+\lambda'_a (a),(\lambda'_a)^2)$. Because $\lambda'_a(Y)=Z$
		and $Z\ne Y$ we have that $\lambda'_a (Z)\ne Z$. Hence $a+\lambda'_a (a)$ is divisible by all
		elements in $Z\cup \lambda'_a (Z)$. As $Z$ is properly contained in $Z\cup \lambda'_a (Z)$ this
		yields $(a,\lambda'_a)^{2}\in M_{i+1}$. Part (1) and Lemma~\ref{lemma1} therefore imply that
		$(M_{*Y}\cup M_{i+1})/M_{i+1}$ is a nil left ideal of $M/M_{i+1}$.
		
		(3) This is proved similarly as part (2).
	\end{proof}
\end{lemma}

\begin{lemma}\label{lemma3}
	Let $Y\in \mathcal{L}$. If $a_Y \notin A_{i+1} $, i.e., the generators of $A$ that divide $a_Y$
	are precisely those that belong to $Y$, then the following properties hold:
	\begin{enumerate}
		\item The derived solution $s\colon\{a_1,\dotsc,a_n\}^{2}\to\{a_1,\dotsc,a_n\}^2$ restricts
		to a solution $s_Y\colon Y^2\to Y^2$ of the YBE.
		\item $M_{YY}$ is a subsemigroup of $M$.
		\item There exists a positive integer $t$, so that for all $k \geq t$, $m_Y^{ k} M_{YY}$ is a cancellative subsemigroup of $M$ and it is an ideal of $M_{YY}$.
		We call it a cancellative component of $M$.
		\item $M_{XX}=M_n$ and $G(X,r)$ is the group of fractions of $m_X^{k} M_{XX}$.
	\end{enumerate}
	In particular, by \eqref{ia_Y}, replacing if necessary $d$ by a multiple, we may assume that $m_YM_{YY}$ is cancellative, for all $Y \in \mathcal{L}$ with $a_Y\notin A_{i+1}$, and $G(X,r)$ is the group of fractions of $m_X M_{XX}$.
	\begin{proof}
		(1) Let $x,y\in Y$ and suppose $s(x,y)=(u,v)$. We need to show that $u,v\in Y$. If $x\ne y$ then
		$a_Y=dx+d y+ b$ with $b=\sum_{x,y\ne z\in Y}dz$. So, \[a_Y=(d-1)x+x+y+(d-1)y+b=(d-1)x+ u+v+(d-1)y+b,\]
		and thus $a_Y $ is divisible by $u$ and $v$. Whereas, if $x=y$ then $a_Y=dx+c$ with $c=\sum_{x\ne z\in Y}dz$.
		So \[a_Y=(d-2)x+x+x+c=(d-2)x+u+v+c,\] and thus $a_Y$ is divisible by $u$ and $v$ as well. Hence the assumption
		$a_Y\notin A_{i+1}$ yields that $u,v\in Y$, as desired. Therefore, (1) follows.
		
		(2) From part (1) it follows that $M_{YY}$ does not contain elements of $A_{i+1}$ and is multiplicatively closed.
		
		(3) From \cite[Proposition~4.2]{JKV} we know that for any solution $(X,r)$ of the YBE, there exists $t \geq 1$, so that for all $k \geq t$, $m_X^{ k} M$ is a cancellative
		ideal of $M$ (equivalently $ka_XA$ is a cancellative ideal of $A$). Because of (1) we know that $(Y,s_Y)$ is a
		solution of the YBE. Hence the structure monoid of $(Y,s_Y)$ has a cancellative ideal generated by $ka_Y$, for some $t\geq 1$ and all $k\geq t$.
		As $M_{YY}$ is a subsemigroup of $M$, we get that, there exists $t\geq 1$ such that for all $k\geq t$, $m_Y^kM_{YY}$ is a cancellative semigroup and 
		an ideal of $M_{YY}$, as desired.
		
		(4) It is clear that $M_{XX}=M_n$. By part (3), there exists a positive integer $t$ such that $S=m_X^{k}M_{XX}$, with $k\geq t$, is a cancellative subsemigroup of $M$. Since $m_X^{k}\in M_{XX}$ and $(a_X+a_i,\lambda'_{a_i})\in M_{XX}$
		for each generator $a_i$ of $A$, we get $m_X^{k+1}\in S$ and $m_X^{k+1}x_i=m_X^k(a_X+a_i,\lambda'_{a_i})\in S$. Hence each
		element $(m_X^{k+1})^{-1} (m_X^{k+1} x_i) $ is in the group of fractions $H=SS^{-1}$ of $S$. Next, observe that the natural
		morphism $S\to G(X,r)$ is injective. Indeed, if $a,b\in S$ are equal in $G(X,r)$ then there exists $l\geq k+1$ such
		that $m_X^{l}  a=m_X^{l}  b$
		in $S$; see \cite[Lemma 6.3]{JKV}. Since $S$ is cancellative and $m_X^l\in S$, we get $a=b$.
		The embedding $S\to G(X,r)$ induces an embedding $H\to G(X,r)$. If $\varphi\colon X\to G(X,r)$ is the natural map
		then, because of the above, $\varphi(x_i)$ is in the image of the embedding $H\to G(X,r)$. As $G(X,r)$ is generated
		by all $\varphi(x_i)$, we conclude that $H=G(X,r)$.
	\end{proof}
\end{lemma}

\begin{lemma}\label{lemma4}
	Let $Y,Z,U,V\in \mathcal{L}$. If $Z\ne U$ then $M_{YZ}M_{UV}\subseteq M_{i+1}$.
	\begin{proof}
		Let $(a,\lambda'_a)\in M_{YZ}$ and $(b,\lambda'_b)\in M_{UV}$. Then $a+\lambda'_a (b)$ is divisible by all
		the elements of $ Y \cup \lambda'_a (U)$. Because $Z\ne U$ and $\lambda_a (Z)=Y$, we have $\lambda'_a (U) \ne Y$.
		Hence, $|Y\cup\lambda'_a (U)|>i$ and thus $(a,\lambda'_a)(b,\lambda'_b)=(a+\lambda'_a (b),\lambda'_a \lambda'_b)\in M_{i+1}$,
		as desired.
	\end{proof}
\end{lemma}

Let \[\mathcal{L}_{u}=\{ Y\in \mathcal{L}:a_Y \notin A_{i+1}\}.\] Because of Lemma~\ref{lemma3},
\[\mathcal{L}_{u}=\{Y\in \mathcal{L}:M_{YY}\text{ is a subsemigroup of }M\}.\]
We define the relation $\sim$ on $\mathcal{L}_u$ as follows. For $Y,Z\in\mathcal{L}_u$ we put
 \[Y\sim Z \text{ if and only if } M_{YZ}\ne \varnothing \text{ or } M_{ZY} \ne \varnothing .\]

\begin{lemma}\label{lemma5}
	The following properties hold:
	\begin{enumerate}
		\item if $Y,Z\in \mathcal{L}_u$ then $Y\sim Z$ if and only if $M_{YZ}\ne \varnothing$ and $M_{ZY} \ne \varnothing$.
		\item $\sim$ is an equivalence relation on $\mathcal{L}_u$.
	\end{enumerate}
	\begin{proof}
		(1) Let $Y,Z\in \mathcal{L}_u$ and $Y\sim Z$. Suppose $M_{YZ}\ne \varnothing$. Let $(a,\lambda'_a)\in M_{YZ}$,
		in particular, $\lambda'_a (Z)=Y$. Since $a_Y\notin A_{i+1}$ there exists $b\in A_{i}\setminus A_{i+1}$ such
		that $(a+\lambda'_a(b),\lambda'_a\lambda'_b)=(a,\lambda'_a) (b,\lambda'_b)=(ka_Y,\id)$ for some positive integer $k$.
		This yields $\lambda'_b=(\lambda'_a)^{-1} $ and thus $\lambda'_b (Y)=Z$. Because of Lemma~\ref{lemma3}, the generators
		of $A$ that divide $\lambda'_a (b)$ are precisely the elements of $Y$. Hence, the generators of $A$ that divide $b$
		are precisely the elements of $(\lambda'_a)^{-1}(Y)=Z$. It follows that $(b,\lambda'_b)=(b,(\lambda'_a)^{-1})\in M_{ZY}$.
		Hence $M_{ZY}\ne \varnothing$. Part (1) then follows.
		
		(2) Clearly $\sim$ is reflexive and symmetric. To show that it is transitive, let $Y,Z,U\in \mathcal{L}_u$ with
		$Y\sim Z$ and $Z\sim U$. So, the sets $M_{YZ},M_{ZY},M_{ZU},M_{UZ}$ are all non-empty. It follows that also all
		of the sets $m_Y M_{YZ}$, $m_ZM_{ZY}$, $m_ZM_{ZU}$, $m_UM_{UZ}$ are non-empty. We need to show that also the set
		$M_{YU}$ is non-empty. To do so, it is sufficient to show that ($m_YM_{YZ})(m_ZM_{ZU})$ is not contained in 
		$M_{i+1}$. Let $(a,\lambda'_a)\in M_{YZ} $ and $(b,\lambda'_b)\in M_{ZU} $. Since $\lambda'_a(Z)=Y$ then
		$\lambda'_a(a_Z)=a_Y$ and thus \[(a_Y+a,\lambda'_a)(a_Z+b,\lambda'_b)=(a_Y+a+a_Y+\lambda'_a(b),\lambda'_a\lambda'_b).\]
		Because of Lemma~\ref{lemma3}, the generators of $A$ that divide $a_Y+a+ a_Y+\lambda'_a (b)$ are precisely the elements
		of $Y$. Furthermore, $\lambda'_a \lambda'_b (U)=\lambda'_a (Z)=Y$. Hence, we get
		$(a_Y+a,\lambda'_a)(a_Z+b,\lambda'_b)\in M_{YU}$. This proves part (2).
	\end{proof}
\end{lemma}

\begin{proposition}\label{prop1}
	Let $\mathcal{L}_1 ,\dotsc,\mathcal{L}_{k}$ denote the equivalence classes of $\sim$ on $\mathcal{L}_u$.
	For each $1\leq j \leq k$ put \[\mathcal{U}_{ij}=\bigcup_{Y,Z\in\mathcal{L}_j}M_{YZ},\quad
	U_{ij}=\bigcup_{Y,Z\in \mathcal{L}_j}m_YM_{YZ},\quad U_i=\bigcup_{j=1}^kU_{ij}.\]
	The following properties hold:
	\begin{enumerate}
		\item $(\mathcal{U}_{ij} \cup M_{i+1})/M_{i+1} $ is a subsemigroup of $M_i/M_{i+1}$ with
		$M_{YZ}M_{ZV}\subseteq M_{YV}$ and $M_{YZ}M_{UV}\subseteq M_{i+1}$ for all $Y,Z,U,V\in \mathcal{L}_{j}$ with $U\ne Z$.
		\item $(U_{ij} \cup M_{i+1})/M_{i+1}$ is an ideal of $(\mathcal{U}_{ij} \cup M_{i+1})/M_{i+1} $ and it is a
		subsemigroup of a completely $0$-simple inverse semigroup with maximal subgroups the group of fractions of $m_YM_{YY}$.
		For simplicity we denote the former as $U_{ij}^{0}$ and we call it a uniform component of $M$ of degree $|\mathcal{L}_j|$.
		\item $(\mathcal{U}_{ij} \cup M_{i+1})/M_{i+1} $ does not contain a nil ideal.
		\item $B_i=M_{i+1}\cup\bigcup_{Y:a_Y\in A_{i+1}}(M_{*Y}\cup M_{Y*})$ is the radical of $M_i/M_{i+1}$.
		Furthermore, if $a_Y\in A_{i+1}$ and $Z\in \mathcal{L}$ with $a_Z\notin A_{i+1}$ then $M_{YZ}=\varnothing$ or $M_{ZY}=\varnothing$.
		\item $M_i/B_i=\bigcup_{j=1}^k(\mathcal{U}_{ij} \cup B_i)/B_i$, a $0$-disjoint union.
		\item $M_i/(U_i \cup M_{i+1})$ is a nil semigroup.
	\end{enumerate}
	Hence we have an ideal chain
	\begin{multline*}
		M_{i+1}\subseteq B_i\subseteq U_{i1}\cup B_i\subseteq U_{i1}\cup U_{i2}\cup B_i\subseteq\dotsb\\
		\dotsb\subseteq U_{i1}\cup U_{i2} \cup \dotsb \cup U_{ik}\cup B_i=U_i\cup B_i\subseteq M_i,
	\end{multline*}
	where the first and last Rees factor is a power nilpotent semigroup and all other Rees factors are uniform
	subsemigroups of an inverse completely $0$-simple semigroup with maximal subgroups the groups of fractions
	of cancellative subsemigroups of $M$.
	\begin{proof}
		(1) Because of Lemma~\ref{lemma3} and Lemma~\ref{lemma4}, we only need to show that $M_{YZ}M_{ZV}\subseteq M_{YV}$.
		To do so, let $(a,\lambda'_a)\in M_{YZ}$ and $(b,\lambda'_b)\in M_{ZV}$. Then
		$(a,\lambda'_a)(b,\lambda'_b)=(a+\lambda'_a(b),\lambda'_a \lambda'_b)$ and $\lambda'_a (b) \subseteq \langle Y\rangle$.
		From Lemma~\ref{lemma3} we get that the generators of $A$ that divide $a+\lambda'_a (b)$ are precisely the elements of $Y$.
		Furthermore, $\lambda'_a \lambda'_b (V)=\lambda'_a (Z)=Y$. Hence $(a,\lambda'_a) (b,\lambda'_b)\in M_{YV}$, as desired.
		
		(2) It is easily verified that $U_{ij}^0=(U_{ij}\cup M_{i+1})/M_{i+1}$ is an ideal of $(\mathcal{U}_{ij}\cup M_{i+1})/M_{i+1}$.
		Because of Lemma~\ref{lemma3}, we know that its diagonal components, i.e., the subsemigroups $m_YM_{YY}$ are cancellative, with
		$Y\in \mathcal{L}_{j}$.
		
		Note that for $a\in A$ we have $(a,\lambda'_a)^{d^2}=(b,\id)^d=(db,\id)$ for some $b\in A$. Hence, for $a,b\in A$,
		the elements of the type $(a,\lambda'_a)^{d^2}$ and $(b,\lambda'_b)^{d^2}$ commute.
		
		It follows that each diagonal component $m_YM_{YY}$ is an Ore semigroup with a group of fractions, denoted $G_{YY}$.
		Actually $G_{YY}$ is obtained from $m_YM_{YY}$ by inverting the element $m_Y$. It is then readily verified that $U_{ij}^{0}$
		is uniform in the inverse completely $0$-simple semigroup $\mathcal{M}^{0}(G_{YY},d_j,d_j,I_j)$, where $d_j=|\mathcal{L}_{j}|$
		and $I_j$ is the identity matrix of degree $d_j$.
		
		(3) Because of Lemma~\ref{lemma3} we know that each $M_{YY}$ is a subsemigroup of $M$, in particular it does not
		contain a nil ideal. Part (3) is now straightforward to check, using standard calculations in the uniform semigroup.
		
		(4) We know from Lemma~\ref{lemma2} that $B_i$ is in the (nil) radical of $M_i/M_{i+1}$. From part (3) we also know
		that the nil radical of $M_i/M_{i+1}$ does not intersect any of the uniform components. Hence $B_i$ is the nil radical
		of $M_i/M_{i+1}$. To prove the second statement of (4). Assume $Z,Y\in \mathcal{L}$ and $a_Y\in A_{i+1}$. Suppose then
		$M_{YZ}\ne \varnothing$ and $M_{ZY}\ne \varnothing$. Then, $M_{ZY}M_{YZ}\subseteq M_{ZZ}$. Since $M_{YZ}\subseteq B_i$
		and because $B_i$ is an ideal, it follows that $M_{ZZ}\subseteq B_i$, a contradiction.
		
		(5) and (6) These are now obvious.
	\end{proof}
\end{proposition}

Let us now verify when $M=M(X,r)$ is Malcev nilpotent.

\begin{theorem}\label{mainnilpotent}
	Let $(X,r)$ be a finite solution of the YBE. Let $n=|X|$. Then, the structure monoid $M=M(X,r)$
	is Malcev nilpotent if and only if all cancellative subsemigroups of $M$ (actually it is sufficient that
	all $m_Y M_{YY}$ with $m_Y\notin M_{|Y|+1}$) are Malcev nilpotent and if furthermore the following condition
	(called the Nilpotency Condition) is not satisfied:
	\begin{gather}
		\begin{gathered}
			\text{there exist subsets $Y\ne Z$ of $\{a_1,\dotsc,a_n\}$, the generators of $A(X,r)$,}\\
			\text{with $a_Y$ and $a_Z$ only divisible by elements of $Y$, respectively $Z$, and}\\
			\text{$a,b\in\free{Y\cap Z}$ such that $\lambda_b'(\lambda_a'^{-1}(Y))=Z$ and $\lambda_b'(\lambda_a'^{-1}(Z))=Y$.}
			\end{gathered}\tag{NC}\label{NC}
		\end{gather}
	\begin{proof}
		Because of the constructed ideal chain of $M$ it follows from \cite[Theorem~11]{JR2006} that $M$ is nilpotent
		if and only if all cancellative components of $M$ are nilpotent and if, furthermore, there does not exist a
		subset $F=\{ f_1,f_1',f_2,f_2'\}$ in a uniform component of $M_i/M_{i+1}$, say $U_{ij}^{0}$, such that all $f_k,f_k'$
		belong to some cancellative component, $f_k$ and $f_k'$ do not belong to the same cancellative component for $k=1,2$,
		and there exist elements $u_1,u_2\in M$ such that the elements $f_2 u_1 f_1$, $f_2'u_2 f_1$, $f_2' u_1 f_1'$ and
		$f_2 u_2f_1'$ are all non-zero in $U_{ij}^0$.
		
		So, in order to prove the result it is sufficient to prove that the existence of such elements can be translated
		into condition \eqref{NC}. Because of Lemma~\ref{lemma3}, without loss of generality, we may assume that each
		$f_1, f_1', f_2, f_2'$ has the permutation coordinate equal to the identity. From Proposition~\ref{prop1} we know
		that there exists distinct subsets $Y$ and $Z$ (of cardinality $i$) of $\{a_1,\dotsc,a_n\}$ with $f_2\in m_YM_{YY}$
		and $f_2'\in m_ZM_{ZZ}$ with $m_Y\notin M_{i+1}$ and $m_Z\notin M_{i+1}$. 
		Also, $f_1\in m_VM_{VV}$ and $f_1'\in m_WM_{WW}$,
		for some distinct subsets $V,W$ of $\{a_1,\dotsc,a_n\}$ (also of cardinality $i$) with $m_V\notin M_{i+1}$ and
		$m_W\notin M_{i+1}$, and $Y,Z,V,W$ are equivalent for the relation $\sim$. Write $u_1=(a,\lambda'_a)$. Because
		$f_2 u_1 f_1\notin M_{i+1}$ and $f_2 u_1 f_1\in M_{YV}$ we get that $a$ can only be divisible by elements
		in $Y$ and $\lambda'_a (V)=Y$. As $f_2' u_1 f_1'\notin M_{i+1}$ and $f_2' u_1 f_1'\in M_{ZW}$ we also get that
		$a$ can only be divisible by elements of $Z$ and $\lambda'_a (W)=Z$. Hence, $a$ can only be divisible by elements
		of $Y\cap Z$ and $a\in \langle Y\cap Z \rangle$. Write $u_2=(b,\lambda'_b)$. The conditions
		$f_2'u_2 f_1\notin M_{i+1}$ and $f_2 u_2 f_1'\notin M_{i+1}$ yield that $b\in \free{Y\cap Z}$, $\lambda'_b(V)=Z$
		and $\lambda'_b (W)=Y$. Hence, condition \eqref{NC} follows. Because of Proposition~\ref{prop1} it easily is
		verified that condition \eqref{NC} implies the existence of $F$ satisfying the required conditions.
	\end{proof}
\end{theorem}

\begin{corollary}\label{cornilpotent}
	Assume that $(X,r)$ is a finite solution of the YBE. If the group $\gr(\lambda_x:x\in X)$
	is of odd order or if the uniform components have degree one then the structure monoid $M=M(X,r)$
	is Malcev nilpotent if and only if all cancellative subsemigroups of $M$ (actually it is sufficient
	that all cancellative components $m_Y M_{YY}$ with $m_Y\notin M_{|Y|+1}$) are Malcev nilpotent.
	\begin{proof}
		If all uniform components are of degree $1$ then each equivalence class $\mathcal{L}_j$ contains
		one element and hence condition \eqref{NC} is trivially not satisfied (no distinct $Y$ and $Z$ exist).
		Hence the result follows in this case. Assume now that the group $\gr(\lambda_x:x\in X)$ is of odd order.
		With notations as in condition \eqref{NC} let $f=\lambda_b'(\lambda_a')^{-1}$. Then $f(Y)=Z$ and $f(Z)=Y$.
		Hence, $f^{2}(Y)=Y$. Since, by assumption, $f$ has odd order, we get that $f(Y)=Y$. Hence $Y=Z$ and
		condition \eqref{NC} is trivially not satisfied.
\end{proof}
\end{corollary}

Example~\ref{NCExample} shows that Corollary~\ref{cornilpotent} does not hold, in general, in case the group
$\gr(\lambda_x : x\in X)$ has even order.

It is easy to give examples of solutions that satisfy condition \eqref{NC}. This can be done via the solutions
$(B,r_B)$ associated to a finite skew left brace $B$ constructed in the following way. Consider the trivial
left braces $A=(\Z/2\Z)^{4}$ and $C=(\Z/2\Z)^{2}$. Let $\alpha\colon C\to\Aut(A)$ be the morphism of groups such that
\begin{align*}
	\alpha(1,0)(a_1,a_2,a_3,a_4) &=(a_2,a_1,a_3,a_4),\\
	\alpha(0,1)(a_1,a_2,a_3,a_4) &=(a_1,a_2,a_4,a_3)
\end{align*}
for all $a_1,a_2,a_3,a_4\in\Z/2\Z$. Let $B=A\rtimes_{\alpha}C$ be the semidirect product of the trivial braces
$A$ and $C$ via $\alpha$. Recall that the addition in $B$ is defined componentwise, i.e.,
\[(a,c)+(a',c')=(a+a',c+c')\] for all $a,a'\in A$ and $c,c'\in C$. Let $e_1,e_2,e_3,e_4$ be the standard basis
of $A$ as a $(\Z/2\Z)$-vector space. Consider the solution $(B,r_B)$ of the YBE associated to the left brace $B$,
and the following subsets of the left derived structure monoid $A(B,r_B)$:
\[Y=\{(e_1,(0,0)),(e_3,(0,0)),(0,(1,0)),(0,(0,1))\}\] and \[Z=\{ (e_2,(0,0)),(e_4,(0,0)),(0,(1,0)),(0,(0,1))\}.\]
Since $A(B,r_B)$ is the free abelian monoid with basis $B$, it is clear that the elements $a_Y$ and $a_Z$ are only
divisible by elements of $Y$, respectively $Z$. Let $a=(0,(1,0))$ and $b=(0,(0,1))$ be the two elements of $Y\cap Z$.
Note that \[\lambda'_b((\lambda'_a)^{-1}(Y))=Z\quad\text{and}\quad\lambda'_b((\lambda'_a)^{-1}(Z))=Y.\]
Hence condition \eqref{NC} is satisfied.

Theorem~\ref{mainnilpotent} easily can be applied on examples. We illustrate this via the following example.
It yields a Malcev nilpotent structure monoid with all cancellative components contained in an abelian group.

\begin{example}[see \cite{SmokVen}]
	Let $X=\{ 1,2,3,4\}$, $\sigma=(1,2)$ and $\tau=(3,4)$. Define $r(x,y)=(\sigma (y),\tau (x))$ for $x,y\in X$.
	Then $(X,r)$ is a solution of the YBE of order $4$. In its structure group we have $1\circ 2=1\circ 1$ and
	$3\circ 4=4\circ 4$. So $G(X,r)=\gr(1,3)$ and the only relation is $1\circ 3=3\circ 1$. Hence $G(X,r)$ is the
	free abelian group of rank $2$; in particular it is nilpotent. Therefore, $m_XM_{XX}$ has a free abelian
	group of rank two as group of fractions. Now,
	\begin{align*}
		A(X,r) &=\free{X\mid x+y=y+\sigma(\tau(x))\text{ for all }x,y\in X}\\
		& =\free{X\mid 1+x=x+2,\,2+x=x+1,\\
		& \phantom{=\free{X\mid}}\,3+x=x+4,\,4+x=x+3\text{ for all }x\in X}.
	\end{align*}
	It is easy to see that
	\begin{align*}
		A(X,r)=\langle 1,2 \mid 1 & +1=1+2=2+2=2+1 \rangle\\
		&+\langle 3,4\mid 3+3=3+4=4+4=4+3\rangle
	\end{align*}
	and we have the extra relations \[1+3=3+2=2+4=4+1\quad\text{and}\quad 1+4=4+2=2+3=3+1.\]
	Notice that all the latter words are in $A_4$. Let $Y=\{1,2\}$ and $Z=\{3,4\}$. Then
	$A_{YY}=A_2\cap\langle 1,2\rangle=1+\langle 1\rangle$ and $A_{ZZ}=A_2\cap \langle 3,4\rangle=3+\langle 3\rangle $;
	both semigroups are cancellative and commutative. Further, $m_YM_{YY}=(a_Y A_{YY})^e=\{(a,\lambda'_a):a\in A_{YY}\}$;
	and $\lambda_a'=\id$ or $\lambda_a'=\sigma$; but $\sigma$ induces also the identity on $A_{YY}$. So $m_YM_{YY}$ is
	abelian and cancellative. Similarly, $m_ZM_{ZZ}$ is abelian and cancellative. Also $A_2\setminus A_3=(1+\free{1})\cup(3+\free{3})$,
	a disjoint union of cancellative semigroups that are orthogonal modulo $A_3$. Further, $A_1\setminus A_2=X$
	(and thus $M_1^2\subseteq M_2$), $A_3=A_4$ and, as said above, $m_XM_{XX}$ is abelian and cancellative.
	Hence all uniform components are of degree $1$ and all cancellative components are abelian. Therefore, it follows
	from Corollary~\ref{cornilpotent} that $M(X,r)$ is a nilpotent semigroup.
\end{example}

The previous example is a solution of Lyubashenko type (see for example \cite{DR1992}), i.e., a finite solution $(X,r)$
with $r(x,y)=(\sigma (y), \tau (x))$ for $x,y\in X$ and some commuting permutations $\sigma$ and $\tau$ on $X$.
We will now deal with all such solutions and determine when they yield nilpotent structure monoids.

\begin{proposition}\label{corlevel1abelian}
	Assume that $(X,r)$ is a Lyubashenko solution with $r(x,y)=(\sigma (y), \tau (x))$ for $x,y\in X$
	and some commuting permutations $\sigma$ and $\tau$ on $X$. Then the structure monoid $M=M(X,r)$
	is nilpotent if and only of $\sigma=c_1^{k_1} \dotsm c_s^{k_s}$ and $\tau=c_1^{1-k_1} \dotsm c_s^{1-k_s}$,
	where $c_1, \dotsc ,c_s$ are disjoint cycles. In this case, all cancellative components are abelian and
	their group of fractions is of rank $1\leq j\leq s$, and all such numbers $j$ can be reached. Furthermore,
	all uniform components have degree~$1$.
	\begin{proof}
		For convenience we will rewrite the solution $r$ as $r(x,y)=(\sigma(y),\sigma^{-1}\gamma (x))$, where
		$\gamma=\sigma \tau\in \Sym(X)$. The disjoint cycle decomposition of $\gamma$ we write as $\gamma=c_1\dotsm c_s$
		and the content of the cycle $c_i$ we denote by $X_i$, a subset of $X$. So $X=X_1\cup \dotsb \cup X_s$,
		a disjoint union. Note that $X_i$ may be a singleton.
		
		We verify when the monoid $M$ is nilpotent, by verifying when the necessary and sufficient conditions of
		Theorem~\ref{mainnilpotent} are satisfied.
		
		Let us first determine the cancellative components of $M$. We begin with the cancellative component determined
		by the set $X$, i.e., $m_XM_{XX}$. Because of Lemma~\ref{lemma3}, its group of fractions is the structure group
		$G=G(X,r)$. Moreover, $G$ also is the structure group of the injectivization of $(X,r)$, i.e.,
		$G=G(\iota(X),r_{\iota(X)})$, where $\iota\colon X\to G$ is the natural mapping and $r_{\iota(X)}$ is the
		restriction of the solution $r_G$ on $G$ to $\iota(X)\times\iota(X)$. Clearly,
		$r(x,\sigma^{-1}(x))=(x,\sigma^{-1}\gamma (x))$. Hence, in $G$, we have $\sigma^{-1} \gamma (x)=\sigma^{-1}(x)$
		and thus in $G$ we have that $\gamma$ is the identity on $\iota(X)$. Therefore, in $G$, we have
		$xy=\sigma(y)\sigma^{-1}(x)$, and thus $(\iota(X), r_{\iota(X)})$ is an involutive solution of the YBE.
		Lemma~\ref{lemmaabelian} yields that if $M$ is nilpotent (and thus also $m_XM_{XX}$ and $G$) then $G$ is abelian.
		The fact that $\gamma=\id$ on $G$ means that if $c_i$ is a cycle of $\gamma$ then all elements in the content
		of $c_i$ are identified in $G$. Moreover, the associated monoid $A(\iota(X),r_{\iota(X)})$ is the free abelian
		monoid on $k$ generators (the number of cycles of $\gamma$). As $G$ is abelian, we need that $\sigma$ is the
		identity when acting on $\iota(X)$. Therefore, on $X$, $\sigma$ must be such that it permutes the contents
		of each $c_i$. So $\sigma (X_i)=X_i$. Now $\sigma$ and $\gamma$ commute. Hence if $c_1=(x_1,\dotsc,x_t)$,
		where $X_1=\{x_1,\dotsc,x_t\}$, then
		\[\gamma=\sigma\gamma\sigma^{-1}=(\sigma (x_1),\dotsc,\sigma(x_t))(\sigma c_2\sigma^{-1})\dotsm(\sigma c_k\sigma^{-1})\]
		and thus $(\sigma(x_1),\dotsc,\sigma(x_t))=c_1$ (similarly for the other cycles). It follows that
		$\sigma=c_1^{k_1}\dotsm c_s^{k_s}$ for some non-negative integers $k_1,\dotsc,k_s$. Thus $\gamma=c_1\dotsm c_s$
		and $\sigma=c_1^{k_1}\dotsm c_s^{k_s}$ and hence $\tau=c_1^{1-k_1}\dotsm c_s^{1-k_s}$.
		
		Now let us look at other possible cancellative components, using the description of the mapping $\sigma$ and $\tau$.
		Because of Lemma~\ref{lemma3} and Lemma~\ref{lemma4}, such a component is determined by a subset $Y$ of $X$, say of
		cardinality $i$, with $m_Y\in M_i\setminus M_{i+1}$. In particular $s_Y$ is a subsolution of $s$. Hence $Y$ must be
		the union of the contents of some cycles of $\gamma$, i.e., the union of some $X_j$. Say $Y=X_{i_1}\cup\dotsb\cup X_{i_l}$.
		Because of the description of $\sigma$ and $\gamma$ this means that $r_Y$ is a subsolution of $r$. Hence, as for the case
		$X$, $m_YM_{YY}$ has a group of fractions $G(Y,r_Y)$ and this must be abelian, and thus also $M(Y,r_Y)$, is nilpotent.
		
		So we have proved that if $M$ is nilpotent then $\gamma=c_1 \dotsm c_s$ and $\sigma=c_1^{k_1}\dotsm c_s^{k_s}$ and thus
		$\tau=c_1^{1-k_1} \dotsm c_s^{1-k_s}$ and conversely for such permutations we have that all cancellative components are abelian.
		
		So it remains to deal with condition \eqref{NC} stated in Theorem~\ref{mainnilpotent}. Let $Y$ be a subset of $X$ as above.
		Then, for any $m=(a,\lambda'_a)\in M$ with $a\in A_{YY}$, we have $\lambda'_a(Y)=Y$. Hence $Y$ is invariant with respect to
		the action of the group $\gr(\lambda_x:x\in X)$. So, $m=(a,\lambda'_a)\in M_{YY}$. In particular, if $Y\sim Z$ (where $\sim$
		is the equivalence relations defined earlier in the section) then $Y=Z$. Hence, all uniform components have degree $1$
		and thus the result follows from Corollary~\ref{cornilpotent}.
	\end{proof}
\end{proposition}

We finish this section with some examples. The first example is a solution $(X,r)$ with abelian structure group $G(X,r)$
while the derived structure group $A_{\gr}(X,r)$ is not nilpotent and it has a uniform component of degree two. However,
the structure monoid $M(X,r)$ is not abelian but it is Malcev nilpotent. 

\begin{example}\label{ab-nonilpotent}
	Let $X=\Z/3\Z$ and define $r(x,y)=(-y,x-y)$ for $x,y\in X$. Then $(X,r)$ is a solution of the YBE.
	The structure group associated to this solution
	\[G=G(X,r)=\gr( 0,1,2 \mid 0\circ 1=2\circ 2=1\circ 0, 0\circ 2=1\circ 1=2\circ 0 )\]
	is an abelian group (note that $0\circ 1=2\circ 2$ implies $2\circ 0\circ 1=2\circ 2\circ 2$
	and thus $1\circ 2\circ 0=1\circ 1\circ 1=2\circ 2\circ 2=2\circ 1\circ 0$ so that $1\circ 2=2\circ 1$).
	The associated derived group is
	\begin{align*}
		A_{\gr} &=A_{\gr}(X,r)=\gr(0,1,2 \mid 0+2=2+1=1+0,\,0+1=1+2=2+0)\\
		& \cong\gr( a,b \mid a+b+a=b+a+b,\, a+b-a=b+a-b,\, 2a+b=b+2a).
	\end{align*}
		Clearly 
	  $$A_{\gr}/\gr(2a,2b)\cong\gr(a,b \mid 2a=2b=3(a+b)=0)\cong S_3$$
	  (note that  $2a=2b$ in $A_{\gr}$, so in fact  $\gr(2a,2b)=\gr(2a)$)
 and thus $A_{\gr}$ is not nilpotent.
	It easily is verified that condition \eqref{NC} is not satisfied and that all cancellative components are
	nilpotent and thus $M=M(X,r)$ is Malcev nilpotent. Note that $M$ has a uniform component of degree two.
	Indeed, let $Y=\{1\}$ and $Z=\{2\}$. Then one has the uniform components consisting of
	\begin{align*}
		M_{YY} &=\{(a,\lambda'_a) :a\in 1+1+\langle 1+1 \rangle \},\\
		M_{ZZ} &=\{(a,\lambda'_a) : a\in 2+2+\langle 2+2 \rangle \},\\
		M_{YZ} &=\{(a,\lambda'_a) :a\in 1+\langle 1+1 \rangle\},\\
		M_{ZY} &=\{(a,\lambda'_a) : a\in 2+\langle 2+2 \rangle \}.
	\end{align*}
\end{example}

In the following example all conditions of Theorem~\ref{mainnilpotent} are satisfied,
so the structure monoid $M(X,r)$ is Malcev nilpotent.

\begin{example}
	Let $X=\Z/4\Z$ and define $r(x,y)=(-y, x+2y)$ for $x,y\in X$. Then $(X,r)$ is a solution of the YBE,
	the structure group is
	\begin{align*}
		G=G(X,r)=\gr(0,1,2,3 \mid 0\circ 1 &=3\circ 2=2\circ 3=1\circ 0,\,0\circ 2=2\circ 0,\\
		0\circ 3 &=1\circ 2=2\circ 1=3\circ 0,\,1\circ 1=3\circ 3),
	\end{align*}
	and the derived structure group is
	\begin{align*}
		A_{\gr}=A_{\gr}(X,r)=\gr(0,1,2,3\mid 0+3 &=3+2=2+1=1+0,\,0+2=2+0,\\
		0+1 &=1+2=2+3=3+0,\,1+3=3+1).
	\end{align*}
	The structure group is nilpotent of class $2$ ($0$ and $2$ are central elements). One easily verifies
	that the conditions mentioned in Theorem~\ref{mainnilpotent} are satisfied. Thus the structure monoid
	$M(X,r)$ is Malcev nilpotent. Note, that for example one easily verifies that there is a degree two
	uniform component built from the sets $Y=\{1\}$ and $Z=\{3\}$.
\end{example}

We now give an example with a structure group that is not nilpotent, as it has $S_3$ as an epimorphic image.
Therefore, not all cancellative components of $M=M(X,r)$ are Malcev nilpotent. Hence the structure monoid
$M$ is not Malcev nilpotent.

\begin{example}\label{example2}
	Let $X=S_3$ and $r(x,y)=(xy^{-1}x^{-1}, xy^2)$ for $x,y\in X$. One can verify that the structure group
	$G(X,r)$ of the solution $(X,r)$ of the YBE has $X=S_3$ as an epimorphic image, and hence it is not nilpotent.
\end{example}

The following is an example with abelian structure group, but not all cancellative components are nilpotent semigroups.
Furthermore, the condition \eqref{NC} holds.

\begin{example}
	Let $X=\{1,2,3,4\}$. Define $\lambda_1=\lambda_2=\rho_1=(3,4)$ and $\lambda_3=\lambda_4=\rho_2=\rho_3=\rho_4=\id$.
	Moreover, let $r(x,y)=(\lambda_x(y),\rho_y(x))$ for $x,y\in X$. Then $(X,r)$ is a solution of the YBE.
	Furthermore, the associated structure monoid
	\begin{align*}
		M=M(X,r)=\free{X \mid 1\circ 2 &=2\circ 1,\,1\circ 3=4\circ 1,\,1\circ 4=3\circ 1,\\
		2\circ 3 &=4\circ 2=2\circ 4=3\circ 2,\,3\circ 4=4\circ 3},
	\end{align*}
	is not abelian. However, the structure group
	\[G(X,r)=\gr(1,2,3 \mid 1\circ 2=2\circ 1,\, 1\circ 3=3\circ 1,\, 2\circ 3=3\circ 2),\]
	is abelian. The derived structure monoid is
	\begin{align*}
		A=A(X,r)=\langle X\mid 1+2 &=2+1,\,1+3=3+1,\,1+4=4+1,\\
		2+4 &=4+2=2+3=3+2,\,3+4=4+3 \rangle.
	\end{align*}
	Since $A(X,r)$ is abelian, we may take $d=2$. Let $Y=\{1,3 \}$ and $Z=\{1,4\}$. Then
	\begin{alignat*}{2}
		a_Y &=1+1+3+3\in A_2 \setminus A_3,\quad & m_Y &=1\circ 1\circ 3\circ 3\in M_2 \setminus M_3,\\
		a_Z &=1+1+4+4\in A_2 \setminus A_3,\quad & m_Z &=1\circ 1\circ 4\circ 4\in M_2 \setminus M_3.
	\end{alignat*}
	We obtain the following non-empty components
	\begin{align*}
		M_{YY} &=\{(a, \lambda'_a) : a\in 1+1+3+\free{1+1,3}\},\\
		M_{YZ} &=\{(a, \lambda'_a) : a\in 1+3+\free{1+1,3}\},\\
		M_{ZY} &=\{(a, \lambda'_a) : a\in 1+4+\free{1+1,4}\},\\
		M_{ZZ} &=\{(a, \lambda'_a) : a\in 1+1+4+\free{1+1,4}\}.
	\end{align*}
	Finally, let $a=1\in\free{Y\cap Z}$ and $b=1+1\in \langle Y \cap Z \rangle$. Then,
	$\lambda'_{1+1}((\lambda'_1)^{-1}(Y))=\lambda'_1(Y)=Z$ and $\lambda'_{1+1}((\lambda'_1)^{-1}(Z))=\lambda'_1(Z)=Y$.
	Hence, condition \eqref{NC} is satisfied. Furthermore, not all cancellative components are Malcev nilpotent.
	Indeed, take $T=\{1,3,4 \}$. Then $M_{TT}=\{(a,\lambda'_a):a\in\free{1,3,4}\}$. If we restrict $r$ to $T \times T$
	we obtain a subsolution $r_T$ with the structure group
	\begin{align*}
		G(T,r_T) &=\gr(T \mid 1\circ 3=4\circ 1,\,1\circ 4=3\circ 1,\,3\circ 4=4\circ 3)\\
		& \cong\gr(3,4 \mid 3\circ 4=4\circ 3) \rtimes \gr(1),
	\end{align*}
	where the action of $1$ is interchanging $3$ and $4$. The group $G(T,r_T)$ is not nilpotent as it contains
	the infinite dihedral group. Hence, $M_{TT}$ is not Malcev nilpotent. So, we found a solution $(X,r)$ where
	$G(X,r)$ is abelian, condition \eqref{NC} is satisfied, but not all cancellative components are nilpotent
	(so the structure monoid $M(X,r)$ is not nilpotent).
\end{example}

We finish this section with an example of a structure monoid for which \eqref{NC} holds and which has abelian
cancellative components. So it is not a nilpotent semigroup.

\begin{example}\label{NCExample}
	Let $X=\{ 1,2,3,4 \}$. Define $\lambda_1=\lambda_2=\rho_1=(3,4)$, $\lambda_3=(2,4),\lambda_4=(2,3)$ and
	$\rho_2=\rho_3=\rho_4=\id$. Moreover, let $r(x,y)=(\lambda_x(y), \rho_y(x))$ for $x,y\in X$. Then $(X,r)$
	is a solution of the YBE. Furthermore, the associated structure monoid
	\begin{align*}
		M=M(X,r)=\free{X \mid 1\circ 2 &=2\circ 1,\, 1\circ 3=4\circ 1,\, 1\circ 4=3\circ 1,\\
		2\circ 3 &=4\circ 2=2\circ 4=3\circ 2=4\circ 3=3\circ 4}
	\end{align*}
	is not abelian. However, the structure group \[G(X,r)\cong\gr(1,2 \mid 1\circ 2=2\circ 1)\]
	is abelian. The derived structure monoid is
	\begin{align*}
		A=A(X,r)=\langle X \mid 1+2 &=2+1,\, 1+3=3+1,\, 1+4=4+1,\\
		2+4 &=4+2=2+3=3+4=4+3=3+2 \rangle.
	\end{align*}
	Since $A$ is abelian, we may put $d=2$. Let $Y=\{1,3 \}$ and $Z=\{1,4\}$. Then
	\begin{alignat*}{2}
		a_Y &=1+1+3+3\in A_2 \setminus A_3,\quad & m_Y &=1\circ 1\circ 3\circ 3\in M_2 \setminus M_3,\\
		a_Z &=1+1+4+4\in A_2 \setminus A_3,\quad & m_Z &=1\circ 1\circ 4\circ 4\in M_2 \setminus M_3.
	\end{alignat*}
	We obtain the following non-empty components
	\begin{align*}
		M_{YY} &=\{ (a, \lambda'_a) : a\in 1+1+3+ \free{1+1,3}\},\\
		M_{YZ} &=\{ (a, \lambda'_a) : a\in 1+3+\free{1+1,3}\},\\
		M_{ZY} &=\{ (a, \lambda'_a) : a\in 1+4+\free{1+1,4}\},\\
		M_{ZZ} &=\{ (a, \lambda'_a) : a\in 1+1+4+\free{1+1,4}\}.
	\end{align*}
	Let $a=1\in\free{Y\cap Z}$ and $b=1+1\in\free{Y \cap Z}$. Then, $\lambda'_{1+1}((\lambda'_1)^{-1}(Y))=\lambda'_1(Y)=Z$
	and $\lambda'_{1+1}((\lambda'_1)^{-1}(Z))=\lambda'_1(Z)=Y$. Hence, condition \eqref{NC} is satisfied. Note that $M_{YY}$
	and $M_{ZZ}$ are abelian. Let us consider all other non-empty subsemigroups $M_{TT}$ with $|T|<4$. If $|T|=1$ then
	these are $\langle a\rangle^e$, with $a\in X=\{1,2,3,4\}$, and clearly $M_{TT}$ is abelian. If $|T|=2$ then the
	only remaining case is $T=\{ 1,2\}$ and $M_{TT}=\langle 1,2\rangle^e$, an abelian semigroup. In case $|T|=3$ there
	is only one such set with $M_{TT} \neq \varnothing$, namely $T=\{ 2,3,4 \}$. Clearly $(T,r|_{T^2})$ is a subsolution
	of $(X,r)$. Hence $M_{TT}$ has an ideal that is cancellative and that has the structure group $G(T,r|_{T^2})$ as its
	group of fractions. It readily is verified that this group is free abelian of rank $1$. Hence, all cancellative
	components of $M$ are abelian and condition \eqref{NC} is satisfied.
\end{example}

\section{Multipermutation Solutions}\label{sec4}

Let $(X,r)$ be a solution of the YBE (recall that it means $(X,r)$ is bijective and non-degenerate).
In this section, we define the retract solution of $(X,r)$ as a generalization of the involutive case
\cite{ESS} and the finite case \cite{LV}. We extend the notion of the multipermutation solution.
Then we go deeper into the study of the relation between a solution $(X,r)$ of the YBE that is a
multipermutation solution and the skew left brace structure of the structure group $G=G(X,r)$ of $(X,r)$,
as well as the relation with the associated solution $(G,r_G)$ of the skew left brace $G$. This link is
extended to the solution $(M,r_M)$ associated to the structure monoid $M=M(X,r)$. In doing so, we will
generalize some results of Gateva-Ivanova and Cameron in \cite{GIC}.

Recall that in \cite{ESS}, Etingof, Schedler and Soloviev introduced the retract relation of involutive
solutions $(X,r)$ of the YBE. This is the binary relation $\sim$ on $X$ defined by $x\sim y$ if and only
if $\lambda_x=\lambda_y$. Then, $\sim$ is an equivalence relation and $r$ induces an involutive solution
$\ov{r}\colon\ov{X}^{2}\to\ov{X}^2$ of the YBE on the quotient set $\ov{X}=X/{\sim}$.

Lebed and Vendramin, in \cite{LV}, generalized this notion to finite solutions $(X,r)$ of the YBE.
In this case, the retract relation is the binary relation $\sim$ on $X$ defined by $x\sim y$ if and only
if $\lambda_x=\lambda_y$ and $\rho_x=\rho_y$. Then $r$ induces a solution $(\ov{X},\ov{r})$ of
the YBE \cite[Lemma 8.4]{LV}. Note that if $r$ is involutive then $\lambda_x=\lambda_y$ if and only if $\rho_x=\rho_y$.

In the following Lemma~\ref{lem:retractissolution} we see that \cite[Lemma 8.4]{LV} is true for solutions $(X,r)$
of the YBE of arbitrary size.

\begin{lemma}\label{lem:retractissolution}
	Let $(X,r)$ be a solution of the YBE. Define an equivalence relation $\sim$ (also denoted $\sim_X$
	to emphasise the set $X$) on $X$ by \[x\sim y\iff \lambda_x=\lambda_y\text{ and }\rho_x=\rho_y.\]
	Then $r$ induces a solution $\ov{r}$ of the YBE on $\ov{X}=X/{\sim}$, by
	\[\ov{r}(\ov{x},\ov{y})=(\ov{\lambda_x(y)},\ov{\rho_y(x)}),\]
	where $\ov{x}$ denotes the $\sim$-class of $x\in X$. This solution is denoted $\Ret(X,r)$
	and is called the retract of $(X,r)$. One says that $(X,r)$ is retractable if $\sim$ is not the trivial relation.
	\begin{proof}
		Let $x,y,z\in X$ be elements such that $x\sim y$. Note that
		\begin{align*}
			\lambda_{\lambda_x(z)}\lambda_{\rho_z(x)} &
			=\lambda_x\lambda_z=\lambda_y\lambda_z=\lambda_{\lambda_y(z)}\lambda_{\rho_z(y)}
			=\lambda_{\lambda_x(z)}\lambda_{\rho_z(y)},\\
			\rho_{\rho_x(z)}\rho_{\lambda_z(x)} &=\rho_x\rho_z=\rho_y\rho_z
			=\rho_{\rho_y(z)}\rho_{\lambda_z(y)}=\rho_{\rho_x(z)}\rho_{\lambda_z(y)},\\
			\lambda_{\lambda_z(x)}\lambda_{\rho_x(z)} &=\lambda_z\lambda_x
			=\lambda_z\lambda_y=\lambda_{\lambda_z(y)}\lambda_{\rho_y(z)}
			=\lambda_{\lambda_z(y)}\lambda_{\rho_x(z)},\\
			\rho_{\rho_z(x)}\rho_{\lambda_x(z)} &=\rho_z\rho_x=\rho_z\rho_y
			=\rho_{\rho_z(y)}\rho_{\lambda_y(z)}=\rho_{\rho_z(y)}\rho_{\lambda_x(z)}.
		\end{align*}
		Hence $\lambda_z(x)\sim\lambda_z(y)$ and $\rho_z(x)\sim\rho_z(y)$. Therefore $\ov{r}$ is well-defined.
		We know that $(X,r^{-1})$ also is a solution of the YBE. We write \[r^{-1}(x,y)=(\hat{\lambda}_x(y),\hat{\rho}_y(x)).\]
		Note that $x\sim y$ if and only if $xy^{-1}\in\Ker(\lambda)\cap\Ker(\rho)\subseteq G(X,r)$. It follows from Remark~\ref{Kernel}
		that $\hat{\lambda}_z(x)\sim\hat{\lambda}_z(y)$ and $\hat{\rho}_z(x)\sim\hat{\rho}_z(y)$, and so the map
		\[\ov{X}^2\to\ov{X}^2\colon(\ov{x},\ov{y})\mapsto (\ov{\hat{\lambda}_x(y)},\ov{\hat{\rho}_y(x)})\]
		is well-defined. Clearly this is the inverse of $\ov{r}$. It is clear that $\ov{r}$ is a bijective set-theoretic
	 	solution of the Yang--Baxter equation. We will prove that it is non-degenerate. Define
		$\lambda_{\ov{x}}\colon\ov{X}\to\ov{X}$ and $\rho_{\ov{x}}\colon \ov{X}\to\ov{X}$
		by $\lambda_{\ov{x}}(\ov{y})=\ov{\lambda_x(y)}$ and $\rho_{\ov{x}}(\ov{y})=\ov{\rho_x(y)}$.
		Clearly, in order to prove that $\lambda_{\ov{x}}$ is bijective, it is enough to prove that it is injective (surjectivity
		of $\lambda_{\ov{x}}$ is an immediate consequence of surjectivity of $\lambda_x$).
	 	
		Note that for every $x,y\in X$, \[\lambda_{\hat{\lambda}_x(y)}(\hat{\rho}_y(x))=x.\] Let $z\in X$ such that
		$y=\hat{\lambda}_x^{-1}(z)$. Then we have \[\lambda_z(\hat{\rho}_{\hat{\lambda}_x^{-1}(z)}(x))=x.\]
		Hence \[\hat{\rho}_{\hat{\lambda}_x^{-1}(z)}(x)=\lambda_z^{-1}(x).\]
		
		Let $x,y,z\in X$ be elements such that $\lambda_z(x) \sim \lambda_z(y)$. We have that
		\[x=\lambda_z^{-1}(\lambda_z(x))=\hat{\rho}_{\hat{\lambda}^{-1}_{\lambda_z(x)}(z)}(\lambda_z(x))
		\sim\hat{\rho}_{\hat{\lambda}^{-1}_{\lambda_z(x)}(z)}(\lambda_z(y))
		=\hat{\rho}_{\hat{\lambda}^{-1}_{\lambda_z(y)}(z)}(\lambda_z(y))=y.\]
		Hence $\lambda_{\ov{x}}$ is injective and thus it is bijective. Similarly one can prove that
		$\rho_{\ov{x}}$ is bijective. Hence $\ov{r}$ is non-degenerate, and the result follows.
	\end{proof}
\end{lemma}

With this equivalence relation $\sim$ at hand, one can now define, as before, multipermutation solutions and their level.

Let $(X,r)$ be a solution of the YBE. Put \[(X_0,r_0)=(X,r) \quad \text{and}\quad (X_n,r_n)=\Ret (X_{n-1},r_{n-1})\]
for $n\geq 1$. Then one says that $(X,r)$ is a multipermutation solution of level $m$, if $|X_m|=1$ and, if $m$ is
positive, $|X_{m-1}|>1$. In this case we write $\mpl(X,r)=m$. In what follows we denote $(X_n,r_n)$ by $\Ret^{n}(X,r)$
for all $n\geq 0$.

\begin{corollary}
	Let $(X,r)$ be a solution of the YBE. If $(X,r)$ is retractable then its left and right derived solutions
	$(X,s)$ and $(X,s')$ are retractable. In particular, if $(X,r)$ is a multipermutation solution of the YBE
	of finite level then so are $(X,s)$ and $(X,s')$.
	\begin{proof}
		Let $x,y\in X$ be two distinct elements such that $x\sim y$. Since $\lambda_x=\lambda_y$ and $\rho_x=\rho_y$,
		by the proof of Lemma~\ref{Permutation}, we have that $\hat{\lambda}_x=\hat{\lambda}_y$ and $\hat{\rho}_x=\hat{\rho}_y$.
		Hence, \[\sigma_x=\lambda_x \hat{\lambda}_x^{-1}=\lambda_y \hat{\lambda}_y^{-1}=\sigma_y\quad\text{and}\quad
		\tau_x=\rho_x \hat{\rho}_x^{-1}=\rho_y \hat{\rho}_y^{-1}=\tau_y,\] which shows that the left derived solution
		$(X,s)$ and the right derived solution $(X,s')$ are retractable.
	\end{proof}
\end{corollary}

Clearly, the reverse implication does not hold. Take $(X,r)$ any irretractable involutive solution of the YBE.
Its derived solution $(X,s)$ is trivial, hence it is a multipermutation solution.

\begin{definition}
    Let $(B,+,\circ)$ be a skew left brace. Define the socle of $B$ as \[\Soc(B)=\{a\in B:a\circ b=a+b=b+a\text{ for all }b\in B\}.\]
    Moreover, we define the socle series for $B$ as follows. Put $\Soc_0(B)=0$ and, for $n\geq 0$, let $\Soc_{n+1}(B)$ denote the unique
    ideal of $B$ containing $\Soc_n(B)$ such that $\Soc_{n+1}(B)/\Soc_n(B)=\Soc(B/\Soc_n(B))$. If there exists a non-negative integer $n$
    such that $\Soc_n(B)=B$ then $B$ is said to have a socle series and the smallest such $n$ is called the socle length of $B$.
\end{definition} 

\begin{remark}\label{rem:bracesolution}
	It is known (see \cite{Ba2018} or \cite[Proposition 1.1.12]{DBThesis}) that, for every skew left brace $B$,
	we have $\Soc(B)=\{a\in B:\lambda_a=\id_B=\rho_a\}$, where $\lambda_a(b)=-a+a\circ b$ and
	\[\rho_a(b)=(\lambda_b(a))^{-1}\circ b\circ a=\lambda_{\lambda_b(a)}^{-1}(-\lambda_b(a)+b+\lambda_b(a))\]
	for all $a,b\in B$. Furthermore, \[r_B\colon B\times B\to B\times B\colon (a,b)\mapsto (\lambda_a(b),\rho_b(a))\]
	is the solution of the YBE associated to the skew left brace $B$ (see \cite[Theorem 3.1]{GV}).
	It is well-known that the maps
	\[\lambda\colon(B,\circ)\to \Aut(B,+)\colon a\mapsto\lambda_a\quad\text{and}\quad\rho\colon(B,\circ)\to \Sym(B)\colon a\mapsto\rho_a\]
	are, respectively, a homomorphism and an anti-homomorphism of groups.
	
	Note that $\Ret^{n}(B,r_B)=(B/\Soc_n(B),r_{B/\Soc_n(B)})$. Hence $(B,r_B)$ is a multipermutation solution
	of the YBE of level $n$ if and only if $B$ has socle length $n$.
\end{remark}

\begin{lemma}\label{morphismretract}
	Assume that $(X,r)$ is a solution of the YBE. If $\iota\colon X\to G(X,r)$ and $\ov{\iota}\colon\Ret(X,r)\to G(\Ret(X,r))$
	are the canonical maps then the rule $\varphi(\iota(x))=\ov{\iota}(\ov{x})$, where $\ov{x}$ denotes the
	equivalence class of $x\in X$ in $\Ret(X,r)$, induces a surjective morphism of solutions $\varphi\colon \Inj(X,r)\to\Ret(X,r)$.
	Moreover, $\varphi$ induces a morphism of groups $\varphi'\colon G(X,r)\to G(\Ret(X,r))$.
	\begin{proof}
		The result clearly holds, if one shows that $\varphi$ is well-defined. As it was shown in \cite{JKV}, for all elements
		$x, y\in X$ satisfying $\iota(x)=\iota(y)$, it holds that $\lambda_x=\lambda_y$. By left-right symmetry, this also 
		shows that $\rho_x=\rho_y$. Hence $\ov{x}=\ov{y}$, which yields
		$\ov{\iota}(\ov{x})=\ov{\iota}(\ov{y})$. As $\varphi\colon\Inj(X,r)\to\Ret(X,r)$ is a morphism
		of solutions and the canonical map $\iota'\colon \Inj(X,r)\to G(\Inj(X,r))=G(X,r)$ is injective, $\varphi$ induces
		a morphism of groups $\varphi'\colon G(X,r)\to G(\Ret(X,r))$.
	\end{proof}
\end{lemma}

Smoktunowicz and Vendramin \cite[Theorem 4.13]{SmokVen} showed that if $(B,r_B)$ is the associated solution
of a finite skew left brace $B$ of size at least two, then the order of $r_B$ is even. The following proposition
shows that this phenomenon also appears for multipermutation solutions.

\begin{proposition}\label{pro:evenordertheorem}
	Let $(X,r)$ be a solution of the YBE with $|X|>1$. If $(X,r)$ is of finite multipermutation
	level and $r$ is of finite order then $r$ is of even order.
	\begin{proof}
		It is sufficient to show the result for multipermutation solutions of level~$1$. Indeed, if $(X,r)$
		is a solution of level $n$ then $\Ret^{n-1}(X,r)$ is a multipermutation solution of level $1$. Furthermore,
		there exists a canonical surjective morphism of solutions $\varphi\colon(X,r)\to\Ret^{n-1}(X,r)$. This entails
		that if $n$ is the order of $r$, the order of the solution $\Ret^{n-1}(X,r)$ is a divisor of $n$. In particular,
		if the latter is even, then so is the order of $r$. Thus assume that $(X,r)$ is of multipermutation level $1$.
		Then, there exists commuting permutations $\sigma$ and $\tau$ of $X$ such that $\lambda_x=\sigma$ and $\rho_x=\tau$
		for each $x\in X$. Suppose, for contradiction's sake, that the order $m=2k+1$ of $r$ is odd. Then,
		$(x,y)=r^m(x,y)=(\sigma^{k+1}\tau^k(y),\sigma^k\tau^{k+1}(x))$. In particular, we obtain that
		$\sigma^{k+1}\tau^k(y)=x$. However, we also have that $(x,x)=r^m(x,x)=(\sigma^{k+1}\tau^k(x),\sigma^{k}\tau^{k+1}(x))$, 
		which shows that $\sigma^{k+1}\tau^k(x)=x=\sigma^{k+1}\tau^k(y)$. Since $\sigma^{k+1}\tau^k$ is a bijection,
		it follows that $x=y$ for all $x,y\in X$, in contradiction with $|X|>1$.
	\end{proof}
\end{proposition}

In view of Proposition~\ref{pro:evenordertheorem}, the question rises whether non-involutive,
injective solutions of finite multipermutation level  exist. The following example illustrates this.

\begin{example}
	Let $X=\{ 1,2,3,4 \}$. Put $\sigma_1=\sigma_2=(3,4)$ and $\sigma_3=\sigma_4=(1,2)$. Then $(X,r)$,
	where $r(x,y)=(\sigma_x(y),x)$ for $x,y\in X$, is a solution of the YBE of multipermutation level $2$.
	Note that the retract of $(X,r)$ is a trivial solution on a set consisting of two elements. Further,
	\begin{align*}
	    G(X,r)=\gr(1,2,3,4\mid1\circ 2 & =2\circ 1,\,3\circ 4=4\circ 3,\,1\circ 3=4\circ 1=2\circ 4=3\circ 2,\\
	    3\circ 1 & =2\circ 3=4\circ 2=1\circ 4).
	\end{align*}
	From the above presentation it
	follows that $2=3\circ 1\circ 3^{-1}$ and $4=1\circ 3\circ 1^{-1}$ in $G(X,r)$ and thus one may check that
	\begin{align*}
	    G(X,r)\cong\gr(1,3\mid 3 & \circ 1\circ 1=1\circ 1\circ 3,\,3\circ 3\circ 1=1\circ 3\circ 3,\\
	    3 & \circ 1\circ 3\circ 1=1\circ 3\circ 1\circ 3).
	\end{align*}
	In particular, $1\circ 1$ and $3\circ 3$ are central
	elements of	$G(X,r)$ and the quotient \[G(X,r)/\gr(1\circ 1,3\circ 3)\cong\gr(a,b\mid a^2=b^2=(ab)^4=1)\cong D_8\]
	is a non-abelian group, and thus also $G(X,r)$ is a non-abelian group. Having this observation in hand
	we are ready to show that the solution $(X,r)$ is injective. First, note that the group $G(X,r)$ admits
	a morphism into the free abelian group on $\{x,y\}$ by mapping both $1,2$ to $x$ and both $3,4$ to $y$.
	It implies that $1\ne 3$, $1\ne 4$, $2\ne 3$ and $2\ne 4$ in $G(X,r)$. Now suppose, on the contrary,
	that $1=2$ in $G(X,r)$. Then also $3=4$ in $G(X,r)$ and thus $G(X,r)\cong\gr(1,3\mid 1\circ 3=3\circ 1)$ would be
	an abelian group, a contradiction. Similarly, on shows that $3\ne 4$ in $G(X,r)$. Therefore, $X$ embeds
	into $G(X,r)$, as claimed. It is also worth to mention that despite $1\ne 2$ in $G(X,r)$ we have
	\begin{align*}
	    1\circ 1\circ 3\circ 3&=1\circ (1\circ 3)\circ 3=1\circ (2\circ 4)\circ 3=(1\circ 2)\circ 4\circ 3
	    \\&=(2\circ 1)\circ 4\circ 3
	    =2\circ (1\circ 4)\circ 3=2\circ (2\circ 3)\circ 3=2\circ 2\circ 3\circ 3,
	\end{align*} which guarantee that $1\circ 1=2\circ 2$ in $G(X,r)$.
\end{example}

The following result is proven in \cite[Lemma 4]{CJO14} for involutive solutions.

\begin{proposition}\label{prop:epilevel}
	Assume that $(X,r)$ and $(Y,s)$ are solutions of the YBE. Then each surjective morphism $\varphi\colon X\to Y$
	of solutions induces a surjective morphism $\ov{\varphi}\colon\Ret(X,r)\to\Ret(Y,s)$ of retracts.
	In particular, if $(X,r)$ is a solution of the YBE of finite multipermutation level $m$ then any homomorphic
	image of $(X,r)$ (for example the injectivization $\Inj(X,r)$) is a solution of finite multipermutation level bounded by $m$.
	\begin{proof}
		Let $x,y\in X$ be such that $x \sim y$, i.e., $\lambda_x=\lambda_y$ and $\rho_x=\rho_y$. Then, for any $z\in X$,
		it follows that \[\lambda_{\varphi(x)}(\varphi(z))=\varphi(\lambda_x(z))=\varphi(\lambda_y(z))=\lambda_{\varphi(y)}(\varphi(z)).\]
		As $\varphi$ is surjective, this implies that $\lambda_{\varphi(x)}=\lambda_{\varphi(y)}$. Analogously, one proves that
		$\rho_{\varphi(x)}=\rho_{\varphi(y)}$. Therefore, the composition $\pi\circ\varphi\colon (X,r)\to\Ret(Y,s)$, where
		$\pi\colon (Y,s)\to\Ret(Y,s)$ is the canonical morphism, induces a surjective morphism of solutions
		$\ov{\varphi}\colon{\Ret(X,r)}\to\Ret(Y,s)$.
	\end{proof}
\end{proposition}

Also subsolutions inherit the property of being multipermutation. Also the following result is proven in \cite[Lemma 5]{CJO14} for involutive 
solutions.

\begin{lemma}\label{lem:multipermsubsolution}
	Let $(X,r)$ be a solution of the YBE and $(Y,s)$ its subsolution, i.e., $Y\subseteq X$, $r(Y^2)\subseteq Y^2$
	and $s=r|_{Y^2}$. If $(X,r)$ is of finite multipermutation level $m$ then $(Y,s)$ is of finite multipermutation
	level bounded by $m$.
	\begin{proof}
		It is clear that $x\sim_Xy$ for some $x,y\in Y$ implies $x\sim_Yy$. In particular, this entails that there
		exists a surjective morphism of solutions $\varphi\colon Y/{\sim_X}\to Y/{\sim_Y}$, where the former is a
		subsolution of $\Ret(X,r)$. Suppose we have shown that the map $\id\colon Y\to Y$ induces a surjective morphism
		of solutions $\varphi_n\colon \ov{Y}_n\to\Ret^{n}(Y,s)$, where $\ov{Y}_n=\{\ov{y}\in\Ret^{n}(X,r):y\in Y\}$.
		Let $x,y\in\ov{Y}_n$ be such that $\ov{x} \sim_{\Ret^n(X,r)} \ov{y}$. As $\varphi_n$ is a morphism of
		solutions, it follows that $\varphi_n(\ov{x}) \sim_{\Ret^n(Y,s)} \varphi_n(\ov{y})$ in $\Ret^{n}(Y,s)$.
		In particular, $\varphi_n$ induces a surjective map $\varphi_{n+1}\colon \ov{Y}_{n+1}\to\Ret^{n+1}(Y,s)$,
		which is by construction a morphism of solutions.
		
		By induction, it follows that if $(X,r)$ is of multipermutation level $m$ then $(Y,s)$ is of multipermutation level at most $m$.
	\end{proof}
\end{lemma}

In case of multipermutation solutions of level $1$ the injectivization turns out to be always involutive.

\begin{proposition}\label{prop:multlevel1involutive}
	Let $(X,r)$ be a solution of the YBE of multipermutation level~$1$.
	Then the injectivization $\Inj(X,r)$ of $(X,r)$ is an involutive solution of the YBE.
	\begin{proof}
		Clearly, the injectivization $\Inj(X,r)$ is still a multipermutation solution of level at most $1$.
		If it is of level $0$ then nothing remains to be shown. Hence, we replace $r$ by its injectivization.
		Write $r(x,y)=(\sigma(y),\tau(x))$ for $x,y\in X$ and some commuting permutations $\sigma$ and $\tau$
		on $X$. Then \[(x,\tau(x))=r(x,\sigma^{-1}(x))\] for each $x\in X$. As the solution $(X,r)$ is injective,
		it follows that $\tau=\sigma^{-1}$, which shows the result.
	\end{proof}
\end{proposition}

In order to show that a solution $(X,r)$ of the YBE is a multipermutation solution if and only if
so is $(M,r_M)$, where $M=M(X,r)$, we need to prove some more intermediate results. First we relate
the retract relations $\sim_X$ and $\sim_M$, introduced in Lemma~\ref{lem:retractissolution}.

\begin{lemma}\label{siminherit}
	Assume that $(X,r)$ is a solution of the YBE. Let $M=M(X,r)$ and $G=G(X,r)$.
	If $x,y\in X$ are such that $x\sim_Xy$ then $x\sim_Gy$ and $x\sim_My$.
	\begin{proof}
		This is easily shown with the argument after the proof of Proposition~\ref{prop:rel}. Indeed,
		we see that $\lambda_x\in\Sym(X)$ determines $\lambda_x\in \Sym(M)$. Similarly one sees that
		$\rho_x\in \Sym(X)$ determines $\rho_x\in \Sym(M)$. Thus $x\sim_X y$ implies that $x\sim_M y$.
		For $G$ the proof is similar: note that $\lambda_x\in \Sym(X)$ induces the map $\lambda^e_x\in\Sym(G)$
		and $\rho_x$ induces the map $\rho^e_x\in\Sym(G)$. Hence $x\sim_X y$ implies that $x\sim_G y$.
	\end{proof}
\end{lemma}

\begin{corollary}\label{cor:multlevelmxr}
	Let $(X,r)$ be a solution of the YBE of finite multipermutation level $m$. Then the solutions
	associated to $M=M(X,r)$ and $G=G(X,r)$ are of finite multipermutation level, bounded by $m+1$.
	\begin{proof}
		First, we may assume that $(X,r)$ is an injective solution.
		Let \[T_n=M(\Ret^n(X,r))\quad\text{and}\quad H_n=G(\Ret^n(X,r))\] for $n\geq 0$.
		In particular, $M=T_0$ and $G=H_0$. We claim that there exist surjective morphisms of solutions
		$\psi_n\colon T_n\to\Ret^n(M,r_M)$ and $\varphi_n\colon H_n\to\Ret^n(G,r_G)$. We shall prove our
		claim by induction on $n$. If $n=1$ then Lemma~\ref{siminherit} implies that if $x, y\in X$ are such
		that $x\sim_X y$ then $ x\sim_M y$ and $x \sim_G y$. Hence, there exist surjective morphisms of
		solutions $\psi_1\colon T_1\to\Ret(M,r_M)$ and $\varphi_1\colon H_1\to\Ret(G,r_G)$, as desired,
		where the latter is well-defined by Lemma~\ref{morphismretract}. Now, suppose that we have
		surjective morphisms of solutions $\psi_n$ and $\varphi_n$ for some $n\geq 1$. Consider $x,y\in\Ret^n(X,r)$
		such that $x\sim_{\Ret^n(X,r)}y$. Then $\psi_n(x) \sim_{\Ret^n(M,r_M)} \psi_n(y)$ and
		$\varphi_n(x)\sim_{\Ret^n(G,r_G)} \varphi_n(y)$, which implies that there exist surjective morphisms
		of solutions $\psi_{n+1}\colon T_{n+1}\to\Ret^{n+1}(M,r_M)$ and $\varphi_{n+1}\colon H_{n+1}\to\Ret^{n+1}(G,r_G)$.
		Hence, the proof of our claim is complete. Since $|{\Ret^m(X,r)}|=1$, we get $T_m\cong\mathbb{N}$ and
		$H_m\cong\mathbb{Z}$. Under this identification, we have surjective morphisms of solutions
		$\psi_m\colon\mathbb{N}\to\Ret^m(M,r_M)$ and $\varphi_m\colon \mathbb{Z}\to\Ret^m(G,r_G)$, where $\mathbb{N}$
		and $\mathbb{Z}$ are considered as trivial solutions. In particular, $|{\Ret^{m+1}(M,r_M)}|=1$ and
		$|{\Ret^{m+1}(G,r_G)}|=1$, showing the result.
	\end{proof}
\end{corollary}

\begin{theorem}\label{prop:charactmultipermmxrgxr}
	Let $(X,r)$ be a solution of the YBE. The following properties are equivalent:
	\begin{enumerate}
		\item The solution $(X,r)$ is of finite multipermutation level.
		\item The associated solution on $M=M(X,r)$ is of finite multipermutation level.
		\item The associated solution on $G=G(X,r)$ is of finite multipermutation level.
	\end{enumerate}
	\begin{proof}
		If $(X,r)$ is a solution of the YBE of finite multipermutation level then, by Corollary~\ref{cor:multlevelmxr},
		the associated solutions on $M$ and $G$ are both of finite multipermutation level.
		
		As $(X,r)$ is a subsolution of $(M,r_M)$, the associated solution on $M$, it follows by
		Lemma~\ref{lem:multipermsubsolution} that if $(M,r_M)$ is of finite multipermutation level
		then $(X,r)$ is of finite multipermutation level.
		
		Lastly, suppose that the associated solution on $G$ is of finite multipermutation level.
		Let $\GG=\GG(X,r)$. The natural map $G\to\GG$ is a surjective homomorphism of solutions from
		$(G,r_G)$ to $(\GG,r_{\GG})$. By Proposition~\ref{prop:epilevel}, $(\GG,r_{\GG})$ is of finite
		multipermutation level. Consider the map $\psi\colon X\to \GG: \; x\mapsto (\lambda_x,\rho_x^{-1})$.
		Then, $\psi$ is a morphism of solutions from $(X,r)$ to $(\GG,r_{\GG})$. Clearly $\psi$ induces an
		injective morphism of solutions $\ov{\psi}\colon \Ret(X,r)\to (\GG,r_{\GG})$. By
		Lemma~\ref{lem:multipermsubsolution}, $\Ret(X,r)$ is of finite multipermutation level and thus
		so is $(X,r)$. This finishes the proof.
	\end{proof}
\end{theorem}

Note that hidden in the proof of Proposition~\ref{prop:charactmultipermmxrgxr}, it was shown that the retraction
of a solution $(X,r)$ of the YBE is a subsolution of the solutions associated to its permutation group
(as a skew left brace). In particular, this entails the following corollary.

\begin{corollary}
	Let $(X,r)$ be a solution of the YBE. If $(X,r)$ is irretractable then $(X,r)$ is an injective solution.
\end{corollary}

\begin{theorem}\label{thm:lengthmultlevel}
	Assume that $(X,r)$ is a multipermutation solution of the YBE of level $m$. Then the group
	$G=G(X,r)$ is solvable of derived length bounded by $m+1$. Moreover, the monoid $A(X,r)$
	is nilpotent of class at most $m+3$ and the group $A_{\gr}(X,r)$ is nilpotent of class at most $m+1$.
	\begin{proof}
		Because of Theorem~\ref{prop:charactmultipermmxrgxr} and Corollary~\ref{cor:multlevelmxr},
		the associated solution $(G,r_G)$ of $G$ is of multipermutation level at most $m+1$.
		Then, by Remark~\ref{rem:bracesolution}, $G$ has a finite socle series of length at most $m+1$.
		As this series also is a subnormal series with abelian factors, the first part of the result follows.
		Moreover, the same series can be considered as a refinement of the upper central series of $A_{\gr}(X,r)$,
		and thus the group $A_{\gr}(X,r)$ is nilpotent of class not exceeding $m+1$.
		
		As before, let $(X,s)$ be the derived solution of $(X,r)$, that is
		\[s(x,y)=(y,\lambda_y(\rho_{\lambda_x^{-1}(y)}(x)))=(y,\sigma_y (x)).\]
		We know that $\GG(X,s)=\gr(\sigma_x:x\in X)$ is an epimorphic image of $A_{\gr}(X,r)$.
		Hence, $\GG(X,s)$ is a nilpotent group of class at most $m+1$. Hence, Proposition~\ref{racknilpotent}
		yields that $A(X,r)$ is nilpotent of class at most $m+3$.
	\end{proof}
\end{theorem}

In \cite[Theorem 6.10]{GIC} Gateva-Ivanova and Cameron proved that if $(X,r)$ is a square-free involutive
solution of the YBE and it is a multipermutation solution of level $m$, then the structure group $G(X,r)$
is solvable of derived length $\leq m$. The following corollary generalizes this result. Moreover, for
square-free solutions it improves the bounds obtained in Theorem~\ref{thm:lengthmultlevel}.

\begin{corollary}\label{Cor:squarefree}
	Let $(X,r)$ be a square-free solution of the YBE. If $(X,r)$ is a multipermutation solution of level
	$m$ then the associated solution $(G,r_G)$ on $G=G(X,r)$ satisfies $m-1\leq\mpl(G,r_G)\leq m$.
	If, furthermore, $(X,r)$ is injective then $\mpl(G,r_G)=m$.
	
	Moreover, the additive group of the skew left brace $G$ is nilpotent of class $\leq m$ and the
	structure group $G$ is solvable of derived length $\leq m$.
	\begin{proof}
		The proof is by induction on $m$. For $m=1$ the solution $r$ is of the form
		$r(x,y)=(\sigma(y),\tau(x))$ for some $\sigma,\tau\in\Sym(X)$ such that
		$\sigma\tau=\tau\sigma$. Since $(X,r)$ is square-free $\sigma(x)=x=\tau(x)$
		for all $x\in X$. Hence $(X,r)$ is the trivial solution. In this case,
		$G$ is a trivial brace and thus the result follows for $m=1$.
		
		Now, let $m>1$ and suppose that the result holds for square-free solutions of multipermutation
		level at most $m-1$. Let $\GG=\GG(X,r)$. We know that the map
		$\Ret(X,r)\to(\GG,r_{\GG})\colon\ov{x}\mapsto(\lambda_x,\rho_x^{-1})$ is an injective morphism
		of solutions. Hence there is a morphism of skew left braces $\varphi\colon G(\Ret(X,r))\to\GG$ such
		that $\varphi(\ov{x})=(\lambda_x,\rho_x^{-1})$ for all $x\in X$. Since $\Ret(X,r)$ is a square-free
		injective solution of multipermutation level $m-1$, it follows by the induction hypothesis that $G(\Ret(X,r))$
		has multipermutation level $m-1$. Since $\varphi$ is clearly surjective, by Proposition~\ref{prop:epilevel},
		$(\GG,r_{\GG})$ is of multipermutation level $\leq m-1$. Since $\Ret(X,r)$ is of multipermutation level $m-1$,
		by Lemma~\ref{lem:multipermsubsolution}, we have that $(\GG,r_{\GG})$ is of multipermutation level $m-1$.
		Furthermore, there are epimorphisms of skew left braces $G\to\GG\colon a\mapsto (\lambda_a,\rho_a^{-1})$ and
		$\GG\to G/\Soc(G)\colon(\lambda_a,\rho_a^{-1})\mapsto \ov{a}$. By Proposition~\ref{prop:epilevel},
		we thus get $m-1\leq\mpl(G,r_G)\leq m$. If, furthermore, $(X,r)$ is injective then, by
		Lemma~\ref{lem:multipermsubsolution}, $(G,r_G)$ is of multipermutation level $m$. Hence the first 
		part of the result follows by induction.
		
		By Remark~\ref{rem:bracesolution}, the second part of the result follows in a similar
		fashion to Theorem~\ref{thm:lengthmultlevel}.
	\end{proof}
\end{corollary}

Actually, one can see that nilpotency gives severe restrictions on the structure of $G(X,r)$.
Moreover, for multipermutation solutions such that $G(X,r)$ is nilpotent, this characterizes the torsion subgroup of $G(X,r)$.

\begin{proposition}\label{prop417}
	Let $(X,r)$ be a finite multipermutation solution of the YBE. If the structure group $G=G(X,r)$
	is nilpotent then the torsion subgroup $T=T(G)$ of $(G,\circ)$ is finite and equal to the additive commutator
	subgroup $[G,G]_+$ of the group $(G,+)=A_{\gr}(X,r)$.
	\begin{proof}
		By Lemma~\ref{lemmaabelian}, $T$ is a finite subgroup of $G$. The inclusion $[G,G]_+\subseteq T$
		is true for any finite solution of the YBE and it was observed in \cite{JKVV}. For completeness'
		sake we repeat the argument. The group $[G,G]_+$ is a characteristic subgroup of $(G,+)=A_{\gr}(X,r)$.
		As $\lambda_x\in \Aut(A_{\gr}(X,r))$ for any $x\in X$, it follows that $[G,G]_+$ is a left ideal, in
		particular a subgroup of $G$. As $A_{\gr}(X,r)$ is a finitely generated finite conjugacy group, it follows
		that $[G,G]_+$ is a finite group. Hence, it is a torsion subgroup of $G$, which shows that $[G,G]_+\subseteq T$.
		
		Suppose that $(X,r)$ is a multipermutation solution of level $m$. We prove that $T=[G,G]_+$ by induction
		on $m$. If $m=1$ then, by Proposition~\ref{prop:multlevel1involutive}, the solution $\Inj(X,r)$ is involutive
		and thus $G\cong G(\Inj (X,r))$ is a torsion-free group and $(G,+)=A_{\gr}(X,r)$ is a free abelian group.
		Hence, the result holds in this case. Assume now that $m>1$ and the result is true for finite solutions 
		of the YBE of multipermutation level $<m$. Consider the natural morphism $\varphi\colon G\to G(\Ret(X,r))$
		of groups defined in Lemma~\ref{morphismretract}. Let $N$ denote its kernel. Clearly, $N \subseteq\Soc(G)$
		and $G/N\cong G(\Ret(X,r))$. Thus $N$ also is a subgroup of $A_{\gr}(X,r)$, i.e., it is an additive subgroup
		of $G$. Then, the induction hypothesis shows that $T \subseteq N\circ[G,G]_+=N+[G,G]_+$. Let $g\in T$.
		There exist $a\in N$ and $b\in [G,G]_+$ such that $g=a\circ b$. Since $[G,G]_+\subseteq T$, we obtain
		$a\in T\cap N$. Because $N\subseteq \Soc(G)$, we have $a^n=na$ for all integers $n$. Thus $\gr(a)=\gr(a)_+$
		is a finite subgroup of $A_{\gr}(X,r)$. Then $(\gr(a)+[G,G]_+)/[G,G]_+$ is a finite subgroup of
		$A_{\gr}(X,r)/[G,G]_+$. By Proposition~\ref{pro:alltorsionaxr}, $A_{\gr}(X,r)/[G,G]_+$ is a torsion-free
		group. Hence, $a\in [G,G]_+$. Therefore, $g\in [G,G]_+$ and, in consequence, $T\subseteq [G,G]_+$.
		Thus $T=[G,G]_+$ and by induction the result follows.
	\end{proof}
\end{proposition}

\begin{corollary}
	Let $(X,r)$ be a finite multipermutation solution of the YBE. If the structure group $G=G(X,r)$
	is nilpotent, then $\ov{G}=G/[G,G]_+$ is a trivial left brace. In particular, the image
	$(\ov{X},\ov{r})$ of $(X,r)$ in $(\ov{G},r_{\ov{G}})$ is a trivial solution.
	\begin{proof}
		By Proposition~\ref{prop417}, it follows that $[G,G]_+$ is a characteristic subgroup of $G$.
		As $[G,G]_+$ is a characteristic subgroup of the additive structure, it follows that $[G,G]_+$
		is an ideal of the skew left brace $G$. Hence, $\ov{G}=G/[G,G]_+$ has a natural skew left brace 
		structure. Furthermore, the additive structure is abelian by construction, thus $\ov{G}$
		is a left brace. Since $\ov{G}$ is a left brace, the solution $(\ov{G},r_{\ov{G}})$
		is involutive. Clearly, there exists a natural epimorphism
		$\varphi\colon G(\ov{X},\ov{r})\to\ov{G}$. Furthermore, there exists a natural
		epimorphism $\psi\colon G(X,r)\to G(\ov{X},\ov{r})$, induced by the morphism of solutions
		$(X,r)\to (\ov{X},\ov{r})$. Note that $[G,G]_+=\gr(x-\sigma_y(x):x,y\in X)$. Recall that
		$\sigma_y(x)$ is used to define the left derived solution of $(X,r)$. Hence, $\psi(x)=\psi(\sigma_y(x))$
		in $G(\ov{X},\ov{r})$ for all $x,y\in X$. Thus, $\psi$ can be factored through an epimorphism
		$\psi_2\colon \ov{G}\to G(\ov{X},\ov{r})$. As both $\varphi\circ\psi_2$ and
		$\psi_2\circ\varphi$ correspond to the identity mapping on the generators of the corresponding groups,
		it follows that both maps are isomorphisms. In particular, $\ov{G}$ can be treated as the structure
		group of $(\ov{X},\ov{r})$. As $G$ is nilpotent, it follows that $\ov{G}$ is nilpotent.
		By \cite[Theorem~2]{CGS}, it follows that $\ov{G}$ is a trivial left brace. In particular,
		this implies that the solution $(\ov{X},\ov{r})$ is trivial.
	\end{proof}
\end{corollary}

Let $(X,r)$ be a finite solution. Let $G=G(X,r)$. As we have seen in the proof of Proposition~\ref{prop417}, 
$[G,G]_+$ is a subgroup of $G$ (see also \cite{JKVV}). A natural question is whether $[G,G]_+$ is the set
$\{g\in G:\text{the multiplicative order of $g$ is finite}\}$. But this is not true as it is shown in the
following example. Crucial is to construct an example of a skew brace $B$ such that its left ideal $[B,B]_{+}$
is not an ideal, i.e., as a multiplicative group it is not a normal subgroup of $(B,\circ)$. Note that $[B,B]_{+}$
is a normal subgroup of $(B,+)$ and thus it is a strong left ideal of the skew left 
brace $B$, in the sense as introduced in \cite{JKVV1}.

\begin{example}
	Consider the trivial brace $A=(\Z/2\Z)^2$. Then the automorphism group $\Aut(A)$ of $A$ is
	isomorphic to the symmetric group of degree $3$. We define $f+g=f\circ g$ for $f,g\in\Aut(A)$.
	Then $(\Aut(A),+,\circ)$ is a skew left brace. Consider the semidirect product 
	$B=\Hol(A)=A\rtimes \Aut(A)$ of skew left braces. Thus
	\begin{align*}
		((a,b),f)+((c,d),g) &=((a+c,b+d), f\circ g),\\
		((a,b),f)\circ ((c,d),g) &=((a,b)+f(c,d), f\circ g)
	\end{align*}
	for $(a,b),(c,d)\in A$ and $f,g\in\Aut(A)$. Then
	$(B,+,\circ)$ is a skew left brace. Note that $[B,B]_+=\{((0,0),\id),((0,0),f),((0,0),f^2)\}$,
	where $f\in\Aut(A)$ is defined as $f(a,b)=(b,a+b)$ for $(a,b)\in A$. Since
	\begin{align*}
		((1,0),\id)^{-1}\circ ((0,0),f)\circ ((1,0),\id) &=((1,0),\id)\circ ((0,1),f)\\
		& =((1,1),f)\notin [B,B]_+,
	\end{align*}
	we have that $[B,B]_+$ is not an ideal of the skew left brace $B$. By \cite[Proposition 3.18]{Ba2018}
	(or \cite[Corollary 2.3.5]{DBThesis}), there exists a finite solution $(X,r)$ of the YBE such that
	$\GG=\GG_{\lambda,\rho}(X,r)\cong B$ as skew left braces. Let $G=G(X,r)$ and let $h_2\colon G\to\GG$
	be the map defined by $h_2(a)=(\lambda_a,\rho_a^{-1})$ for all $a\in G$. We know (cf. Lemma~\ref{Permutation}
	and Remark~\ref{rem:h2}) that $h_2$ is an epimorphism of skew left braces and $\Ker(h_2)$ is an ideal of
	the skew left brace $G$ contained in its socle. Note that $h_2^{-1}([\GG,\GG]_+)=[G,G]_++\Ker(h_2)$. Since
	$[\GG,\GG]_+\cong[B,B]_{+}$ is not an ideal of the skew left brace $\GG$, we have that $[G,G]_++\Ker(h_2)$
	is not an ideal of the skew left brace $G$. 
	
	Note that if \[[G,G]_+=\{g\in G:\text{the multiplicative order of $g$ is finite}\},\] then $[G,G]_+$ is an
	ideal of the skew left brace $G$ and hence $[G,G]_++\Ker(h_2)$ also is an ideal of $G$, a contradiction.
	Therefore, $[G,G]_+$ is different from \[\{g\in G:\text{the multiplicative order of $g$ is finite}\}.\]
\end{example}

\subsection*{Acknowledgments}

\noindent The first author was partially supported by the MINECO-FEDER MTM2017-83487-P and AGAUR 2017SGR1725 (Spain) grants.
The second author is supported in part by Onderzoeksraad of Vrije Universiteit Brussel and Fonds voor Wetenschappelijk Onderzoek
(Belgium). 
The third author is supported by National Science Centre grant 2020/39/D/ST1/01852 (Poland).
The fourth author is supported by Fonds voor Wetenschappelijk Onderzoek (Flanders), via an FWO post-doctoral fellowship.
The fifth author is supported by Fonds voor Wetenschappelijk Onderzoek (Belgium), via an FWO Aspirant-mandate.

\bibliographystyle{amsplain}
\bibliography{refs}

\end{document}